\documentclass{article}
\usepackage{graphicx} % Required for inserting images
\usepackage{amsthm} % for theorem environments
\usepackage{amsmath}
\usepackage{amsfonts}
\usepackage{amssymb}
\usepackage{mathtools}
\usepackage{scalerel}
\usepackage{physics}
\usepackage{bm}
\usepackage{changepage}
\usepackage{dsfont}
\usepackage{mathrsfs}
\usepackage{float}
\usepackage{listings}
\graphicspath{ {./images/} }
\usepackage{hyperref}

\title{The EM Algorithm in Information Geometry}
\author{Sammy Suliman}
\date{February 2024}

\begin{document}

\maketitle

\begin{abstract}
The purpose of this thesis is to convey the basic concepts of information geometry and its applications to non-specialists and those in applied fields, assuming only a first-year undergraduate background in calculus, linear algebra, and probability theory / statistics. We first begin with an introduction to the EM algorithm, providing a typical use case in Python, before moving to an overview of basic Riemannian geometry. We then introduce the core concepts of information geometry and the $em$ algorithm, with an explicit calculation of both the $e$ and $m$ projection, before closing with a discussion of an important application of this research to the field of deep learning, providing a novel implementation in Python.
\end{abstract}

\section{The EM algorithm}
\textbf{Definition 1.1:} Probability distribution \\
A probability distribution is a non-negative function that either sums or integrates to 1 used to describe the chances of a particular event occurring in the sample space of possible outcomes. A continuous function is known as a \textit{probability density function} (pdf). A discrete function is known as a \textit{probability mass function} (pmf). \\\\
In common applications of mathematics such as machine learning, a common problem involves finding the model (distribution with certain parameters) that best fits the available data for the purposes of inference and prediction. \\\\
\textbf{Definition 1.2:} Likelihood function ($\mathcal{L}$) \\
The probability mass or density function of some observed data, viewed as a function of its unknown parameters. \\\\
Essentially, given data, we can choose the log-likelihood ($l$) function that best models the distribution, then optimize the parameters by setting the first-order partial derivatives with respect to each parameter to 0 and solving. This procedure is known as \textbf{Maximum Likelihood Estimation}. \\\\
For an example, we take the 1-dimensional normal distribution:
$$p(x|\theta) = \frac{1}{\sigma\sqrt{2\pi}} e^{\frac{-(x-\mu)^2}{2\sigma^2}}$$
For a sample consisting of $n$ independent, identically-distributed (i.i.d.) random variables:
$$\mathcal{L}(\mu, \sigma) = \frac{1}{2\pi\sigma^2}^{\frac{n}{2}} \exp{\frac{-\sum^n_{i=1} (x_i - \mu)^2}{2\sigma^2}}$$
$$l(\mu, \sigma) = \frac{-n}{2}\log(2\pi\sigma^2) - \frac{1}{2\sigma^2}\sum^n_{i=1}(x_i - \mu)^2$$
$$0 = \frac{\partial l}{\partial \mu} = \frac{2n(\bar{x} - \mu)}{2\sigma^2}$$
$$\hat{\mu} = \bar{x}$$
$$0 = \frac{\partial l}{\partial \sigma} = \frac{-n}{\sigma} + \frac{1}{\sigma^3}\sum^n_{i=1}(x_i - \mu)^2$$
$$\hat{\sigma}^2 = \frac{1}{n}\sum^n_{i=1}(x_i - \bar{x})^2$$
Maximum-likelihood estimation works well in most basic situations, but there are 2 situations extremely common in real-world application in which it fails, when there is missing data and when there is a mixture of distributions. \\\\
The \textbf{Expectation-Maximization (EM) algorithm} (Depster, Laird, Rubin, 1977) is an iterative procedure used to find the maximum likelihood estimate of a given model in these cases. \\
Given an observed \textbf{y}, and a natural choice for a distribution $f(x | \phi)$, where $\phi$ are the unknown parameters for our model, and the incomplete data specification $g(y | \phi) = \int_{\mathcal{L(y)}} f(x | \phi) \, dx$, the algorithm is given in 2 steps: \\\\
\textbf{Expectation step} \\
Find the statistic that contains the maximum possible information about the parameters ($p$ here represents the p-th step of the algorithm): \\\\
\indent $Q(t \mid t^{(p)}) = \mathbb{E}_{f, t^{(p)}} [\log L(t; X, f)]$ \\\\
\textbf{Maximization step} \\
Find the solution $\phi^{(p+1)}$ such that: \\\\
\indent $t^{(p+1)} = \arg\max_{\theta} Q(t \mid t^{(p)})$ \\\\
Starting from some initial guess for the parameters, the E-step determines which classes the missing data was most likely to fall into by finding the expectation of the log-likelihood of our observed data, then the M-step applies maximum likelihood estimation on the "filled-in" data. \\\\
In the case where we have data arising from a mixture of different distributions, we want to fill in not missing data points, but the unknown, or latent, variables underlying our distributions. The E-step then computes the posterior probability (responsibility) of each latent variable for generating each data point. It calculates how likely it is that each data point was generated by each component given the current parameter estimates. The M-step once again maximizes the likelihood given the posteriors. \\\\
We will do an example on the space of Gaussian mixture models. \\\\
\textbf{Definition 1.3:} General Mixture Model (GMM) \\
A linear combination of N random variables that are observed, with K weights (summing to 1), where the components belong to the same distribution, but with different parameters. \\ \\
For an example, we are in a classroom of math and statistics majors, and we want to look at their GPAs. Both populations have normally-distributed GPAs. So we have a set of normal random variables, in one of two classes (fixed distributions), representing different subpopulations. The distributions of the classes $c_i$ is Bernoulli.
$$P(X = x) = p^{x | c_i=M} (1 - p)^{1 - x | c_i=S}$$
The conditional distributions in each class can be written:
$$P(Y_i = y_i \mid C_i = c) = \mathcal{N}(y_i ; \mu_c, \sigma^2) \cdot {1}(C_i = c)$$
We can think of the class labels as being a set of unknown variables $C = \{ C_M, C_S \}$, the GPAs of the students in the classroom as being the observed data $Y=\{ Y_1, ..., Y_n \}$, and the parameters $\theta$ we are interested in finding are the mean grades of the students in both classes, $\mu_M, \mu_S$. \\ \\
\textbf{Theorem 1.1} Bayes' law \\
The probability of event B given that A has already occurred can be given by:
$$\mathbb{P}(B|A) = \frac{\mathbb{P}(A,B)}{\mathbb{P}(A)} = \frac{\mathbb{P}(B)\mathbb{P}(A|B)}{\mathbb{P}(A)}$$
Therefore:
$$p(C_i = c \mid Y_i = y_i) = \frac{\prod p_{c_i} \cdot \mathcal{N}(y_i; \mu_M, \sigma^2)^{{1}(c = c_i)}}{\sum_{c \in \{M,S\}} p_{c_i} \cdot \mathcal{N}(y_i; \mu_c, \sigma^2)} = q(C_i = c)$$ \\
We can solve for the desired parameters (i.e, $\mu_M, \sigma_M$) given the posterior probabilities. \\\\
This should shed light on the motivation behind the EM algorithm. We can only compute the posterior probabilities given the parameters, and vice versa. When we fix estimates for the parameters and solve for the posteriors, this is the E-step. When we fix estimates for the posteriors and solve for the parameters, this is the M-step. We iterate until the estimates stop changing. \\\\
\textbf{Theorem 1.2 (Jensen's inequality):} $\log(\mathbb{E}[X]) \geq \mathbb{E}[\log(X)]$ \\\\
We want to maximize the log-likelihood $\log p_Y(y; \theta)$, split the probability over each class and multiply by a factor of $\frac{q_C(c)}{q_C(c))}$ where $q_C(c_i)$ is the posterior probability that the class is $c_i$ (expression above), so:
$$\log p_Y(y; \theta) = \log \left( \frac{\sum_{c} p_{Y,C}(y, c; \theta) \cdot q_C(c)}{q_C(c)} \right) = \log \left( \mathbb{E}_{C} \left[ \frac{p_{Y,C}(y, C; \theta)}{q_C(C)} \right] \right)$$
Use Jensen's inequality: $$ \geq \mathbb{E}_{C} \left[ \log \left( \frac{p_{Y,C}(y, C; \theta)}{q_C(C)} \right) \right]$$
Applying the property of joint probabilities; $P(A, B) = P(A) \cdot P(B \mid A)$
\begin{align*}
&= \mathbb{E}_{C} \left[ \log \frac{p_{Y,C}(y, C; \theta)}{q_C(C)} \right] \\
&= \mathbb{E}_{C} \left[ \log \frac{p_{Y}(y; \theta) \cdot p_{C|Y}(C \mid y; \theta)}{q_C(C)} \right] \\
&= \log p_Y(y; \theta) - \mathbb{E}_{C} \left[ \log \frac{q_C(C)}{p_{C|Y}(C \mid y; \theta)} \right] \\
&= g(\theta)
\end{align*}
Maximizing with respect to $q_C$ (posteriors) gives the E-step. We can rearrange the above expression as:
$$\mathbb{E}_{C} \left[ \log p_{Y,C}(y, C; \theta) \right] - \mathbb{E}_{C} \left[ \log q_{C}(C) \right]$$
Maximizing with respect to $\theta$ gives the M-step. Any value of $\theta$ that increases $g(\theta)$ beyond $g(\theta_{old})$ must also increase $l(\theta; X )$ beyond $l(\theta_{old}; X )$. The M-step finds such a $\theta$ by maximizing $\mathbb{E}_{C} \left[ \log p_{Y,C}(y, C; \theta) \right]$ over $\theta$ which is equivalent to maximizing $g(\theta | \theta_{old})$ over $\theta$.\\ \\
Let's compute this out manually in the case of our Gaussian Mixture model: \\
The likelihood: \\
$$L(\theta; X, Y) = \prod_{i=1}^{2} \frac{p_{Y_i}}{\sqrt{2\pi}\sigma_{Y_i}} \exp \left( -\frac{(x_i - \mu_{Y_i})^2}{2\sigma_{Y_i}^2} \right)$$
E-step: \\
$$\mathbb{E} \left[ l(\theta; X, Y) \mid X, \theta_{\text{old}} \right] = \sum_{i=1}^{n} \sum_{j=1}^{m} P(Y_i = j \mid x_i, \theta_{\text{old}}) \ln f_{X|Y}(x_i \mid Y_i = j, \theta) P(Y_i = j \mid \theta)$$
Let's rewrite $P(Y_i = j \mid x_i, \theta_{\text{old}})$ in terms of the posterior probabilities: \\
$$P(Y_i = j \mid x_i, \theta_{\text{old}}) = \frac{f_{X|Y}(x_i \mid Y_i = j, \theta_{\text{old}}) P(Y_i = j \mid \theta_{\text{old}})}{\sum_{k=1}^{m} f_{X|Y}(x_i \mid Y_i = k, \theta_{\text{old}}) P(Y_i = k \mid \theta_{\text{old}})}$$
So we simply compute the following analytically:

$$\mathbb{E} \left[ l(\theta; X, Y) \mid X, \theta_{\text{old}} \right] =
$$
$$
\sum_{i=1}^{n} \sum_{j=1}^{m} \frac{f_{X|Y}(x_i \mid Y_i = j, \theta_{\text{old}}) P(Y_i = j \mid \theta_{\text{old}})}{\sum_{k=1}^{m} f_{X|Y}(x_i \mid Y_i = k, \theta_{\text{old}}) P(Y_i = k \mid \theta_{\text{old}})} \ln f_{X|Y}(x_i \mid Y_i = j, \theta) P(Y_i = j \mid \theta)$$

M-step: \\\\
We maximize $\mathbb{E} \left[ l(\theta; X, Y) \mid X, \theta_{\text{old}} \right]$ by setting the partial derivatives equal to zero and solving for the parameters, $\mu_j, \theta_j^2$: \\

$$\mathbb{E} \left[ l(\theta; X, Y) \mid X, \theta_{\text{old}} \right] = \sum_{i=1}^N \sum_{j=1}^M P(Y_i = j | X_i, \theta_{\text{old}}) \cdot \ln \left[ \frac{1}{\sqrt{2\pi \sigma_j^2}} e^{\frac{-(x_i - \mu_j)^2}{2\sigma_j^2}} P(Y_i=j | \theta) \right]$$
\begin{adjustwidth}{-1cm}{0cm}
$$0 = \frac{\partial Q}{\partial \mu_j} = \sum_{i=1}^N \sum_{j=1}^M P(Y_i = j | X_i, \theta_{\text{old}}) \cdot \frac{1}{\frac{1}{\sqrt{2\pi \sigma_j^2}} e^{\frac{-(x_i - \mu_j)^2}{2\sigma_j^2}} P(Y_i=j | \theta)} \cdot \frac{1}{\sqrt{2\pi \sigma_j^2}} e^{\frac{-(x_i - \mu_j)^2}{2\sigma_j^2}} \cdot \frac{x_i-\mu_j}{\sigma_j^2} \cdot P(Y_i=j | \theta)$$
\end{adjustwidth}
$$0 = \sum_{i=1}^N P(Y_i = j | X_i, \theta_{\text{old}}) \cdot (x_i - \mu_j)$$
$$\sum_{i=1}^N P(Y_i = j | X_i, \theta_{\text{old}}) \cdot \mu_j = \sum_{i=1}^N P(Y_i = j | X_i, \theta_{\text{old}}) \cdot x_i$$
$$\mu_j = \frac{\sum_{i=1}^N P(Y_i = j | X_i, \theta_{\text{old}}) \cdot x_i}{\sum_{i=1}^N P(Y_i = j | X_i, \theta_{\text{old}})}$$
\begin{adjustwidth}{-1cm}{0cm}
$$0 = \frac{\partial Q}{\partial \sigma_j} = \sum_{j=1}^M P(Y_i = j | X_i, \theta_{\text{old}}) \cdot \frac{1}{\frac{1}{\sqrt{2\pi \sigma_j^2}} e^{\frac{-(x_i - \mu_j)^2}{2\sigma_j^2}} P(Y_i=j | \theta)} \cdot \frac{1}{\sqrt{2\pi \sigma_j^2}} e^{\frac{-(x_i - \mu_j)^2}{2\sigma_j^2}} \cdot \frac{(x_i - \mu_j)^2}{\sigma_j^3} \cdot P(Y_i=j | \theta)$$
\end{adjustwidth}
$$\sigma_{j}^2 = \frac{\sum_{i=1}^{n} (x_i - \mu_j)^2 P(Y_i = j \mid x_i, \theta_{\text{old}})}{\sum_{i=1}^{n} P(Y_i = j \mid x_i, \theta_{\text{old}})}$$
See Appendix A for a practical demonstration of how to use the EM algorithm for this example. \\\\
We introduce a statistical distance commonly used to find how 2 probability distributions differ from one another. \\\\
\textbf{Definition 1.4:} Kullback-Leibler divergence\\
$D_{\text{KL}}(P \| Q) = \sum_x P(x) \log \left( \frac{P(x)}{Q(x)} \right)$ (for discrete distributions) \\
$D_{\text{KL}}(P \| Q) = \int P(x) \log \left( \frac{P(x)}{Q(x)} \right)$ (for continuous distributions) \\\\
Notice that the expression of the E-step can be rewritten as: \\
$$\log p_Y(y; \theta) - D_{\text{KL}}(q_C(\cdot) \| p_{C|Y}(\cdot \mid y; \theta))$$
This will have enormous consequences when we discuss the geometric interpretation of the EM algorithm. First, we will need to develop a basic grounding in differential geometry. \\
\section{Introduction to Riemannian Geometry}
\textbf{Definition 2.1:}  manifold \\
A space that locally resembles Euclidean space near each point (is homeomorphic to an open subset of Euclidean space around a neighborhood of each point) \\\\
Riemannian geometry is the study of 'smooth' (i.e, differentiable) manifolds equipped with an inner product operation (metric). For centuries, mathematicians and scientists constrained their vision of geometry (from the Ancient Greek \textit{geo} (earth) \textit{metron} (to measure)) to the three-dimensional Euclidean space that we reside in. But as in many other branches of mathematics, by abstracting and generalizing intuitive concepts, we can unlock fertile areas of research and uncover hidden truths about our universe. For example, Albert Einstein earned his place in history by showing that our universe is only locally Euclidean (i.e, a manifold). His study of the geometric properties of this manifold is now known as the theory of general relativity. For our purposes, we will be looking at the ability to describe concepts from statistics in the language of geometry. \\\\
But how is it possible to describe concepts as theoretical as those in probability theory geometrically? While we may think of geometry as describing physical objects in 3 dimensions that we can see, like the length of a line segment or the volumes of a shape, we can extend this to any vector space on which we can define an inner product operation. Here is a brief reminder of the standard Euclidean inner (dot) product we are already familiar with: \\\\
\textbf{Definition 2.2:} dot product \\
Given 2 vectors, $\mathbf{x}, \mathbf{y}$ of dimension $n$: \\
$\mathbf{x} \cdot \mathbf{y} = \sum_{i=1}^n x_iy_i$ \\
Alternatively, $\mathbf{x} \cdot \mathbf{y} = |x| \cdot |y| \cdot \cos{\theta}$, where $\theta$ is the angle between vectors $\mathbf{x}, \mathbf{y}$, $|x| = \sqrt{\mathbf{x} \cdot \mathbf{x}}$ \\\\
In general: \\\\
\textbf{Definition 2.3:} Hilbert space \\
A vector space equipped with an inner product operation $< \cdot, \cdot >$ that takes in 2 vectors and is conjugate symmetric, linear in the first argument, and positive-definite, and that is analytically complete. \\\\ 
Therefore, we can define the notion of an angle in an arbitrary Hilbert space as below:
$$\theta = \arccos(\frac{<\mathbf{x}, \mathbf{y}>}{|x| \cdot |y|})$$
and similarly for distance:
$$d(\mathbf{x}, \mathbf{y}) = |\mathbf{x} - \mathbf{y}|$$
We now can extend the ideas of geometry to any general Hilbert space. \\\\
When it comes to manifolds, one geometric property that is of great importance is its curvature. For a surface (embedded in 3-dimensional Euclidean space), the curvature can change depending on which curve along the surface it is measured. In particular, we are concerned with the minimal and maximal curvatures of a surface. Their product is known as \textbf{Gaussian curvature} and is a fundamental tool in studying surfaces. \\\\
We will do an example computation on the object known as the torus. \\\\
\textbf{Definition 2.4:} torus \\
\begin{align*}
x &= (R + r \cdot \cos \theta) \cdot \cos \phi \\
y &= (R + r \cdot \cos \theta) \cdot \sin \phi \\
z &= r \cdot \sin \theta
\end{align*}
where $R$ is the outer radius, $r$ is the inner radius. \\\\
\begin{figure}[htb]
\centering
\includegraphics[width=5cm, height=3cm]{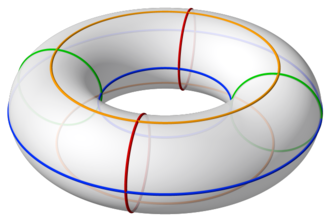}
    \begin{center}
        Figure 2.1: The torus parameterized in $\mathbb{R}^3$. Credit: Wikipedia
    \end{center}
\end{figure}
\\
If we want to understand the curvature of the curves along the torus, we have to start by computing the tangent vectors to the curves, which lie on the tangent plane. A basis for the tangent plane (at any point on the torus), is given by $\{ x_\theta, x_\phi \}$.
$$x_\theta = (-r\cos{\phi}\sin{\theta}, -r\sin{\phi}\sin{\theta}, r\cos{\theta})$$
$$x_\phi = (-(R + r\cos{\theta})\sin{\phi}, (R + r\cos{\theta})\cos{\phi}, 0)$$
Remembering that the dot product is the metric element on a manifold such as the tangent plane, we compute the following matrix:
$$\begin{pmatrix}
E & F \\
F & G \\
\end{pmatrix}
=
\begin{pmatrix}
\mathbf{x}_\phi \cdot \mathbf{x}_\phi & \mathbf{x}_\phi \cdot \mathbf{x}_\theta \\
\mathbf{x}_\phi \cdot \mathbf{x}_\theta & \mathbf{x}_\theta \cdot \mathbf{x}_\theta \\
\end{pmatrix}
=
\begin{pmatrix}
(R + r\cos^2{\theta}) & 0 \\
0 & r^2 \\
\end{pmatrix}
$$
This is known as the \textbf{first fundamental form} matrix. \\\\
Intuitively, if we wanted to understand the curvature of a three-dimensional surface embedded in Euclidean space, what would you want to know? You would want to see how the surface varies across space. This can be quantified as how the normal vector ($\mathbf{n} = \mathbf{x}_u \cross \mathbf{x}_v$)to the surface changes in the direction of each basis vector for the tangent plane. This provides some motivation for the \textbf{second fundamental form} matrix:
$$\begin{pmatrix}
L & M \\
M & N \\
\end{pmatrix}
=
\begin{pmatrix}
\mathbf{x}_{\phi\phi} \cdot \mathbf{n} & \mathbf{x}_{\phi\theta} \cdot \mathbf{n} \\
\mathbf{x}_{\phi\theta} \cdot \mathbf{n} & \mathbf{x}_{\theta\theta} \cdot \mathbf{n} \\
\end{pmatrix}
=
\begin{pmatrix}
\frac{-\cos{\theta}}{r} & 0 \\
0 & \frac{-1}{R + r\cos{\theta}} \\
\end{pmatrix}
$$
It is clear that there are 2 eigenvalues for this matrix, $\lambda_1 = \frac{-\cos{\theta}}{r}, \lambda_2 = \frac{-1}{R + r\cos{\theta}}$, and that these eigenvalues therefore represent the minimal and maximal curvatures of the torus. Together, they are known as the \textit{principal curvatures}. \\\\
Therefore, the Gaussian curvature for a torus is:
$$K = \frac{\cos{\theta}}{r(R + r\cos{\theta})}$$
Depending on the value of $\theta$, the Gaussian curvature may be negative (on the interior of the torus), positive (on the exterior of the torus), or zero (on the top and bottom of the torus, corresponding to $\theta$, the angle of rotation about the tube, equalling $\frac{\pi}{2}$). \\\\
A theorem of fundamental importance in differential geometry which we would be remiss not to mention is the following: \\\\
\textbf{Theorem 2.1:} Theorema Egregium \\
The Gaussian curvature of a surface is invariant under a local isometry (distance-preserving transformation). \\\\
This means that if a surface is bent without stretching or otherwise altering the local angles and distances on the surface, its Gaussian curvature will be preserved. Another way of stating this is that Gaussian curvature is independent of a surface's embedding in ambient space, or is an \textit{intrinsic invariant} of the surface. One application of this is that the Gaussian curvature of a cylinder can automatically be determined to be zero, since a cylinder can be unfolded into a flat plane without stretching, which trivially has curvature zero (the formal computation of this is best left as an exercise to the reader). This is remarkable because the way we defined Gaussian curvature was as the product of the eigenvalues of the second fundamental form matrix, which is an extrinsic measure on a surface, since it involves the normal vector, which lies within the ambient space. But in fact, the Gaussian curvature can be computed completely in terms of the first fundamental form: \\
$$K = \frac{
  \begin{vmatrix}
    -\frac{1}{2}E_{\theta\theta} + F_{\phi\theta} - \frac{1}{2}G_{\phi\phi} & \frac{1}{2}E_{\phi} & F_{\phi} - \frac{1}{2}E_{\theta} \\
    F_{v} - \frac{1}{2}G_{\phi} & E & F \\
    \frac{1}{2}G_{\theta} & F & G
  \end{vmatrix}
  -
  \begin{vmatrix}
    0 & \frac{1}{2}E_{\theta} & \frac{1}{2}G_{\phi} \\
    \frac{1}{2}E_{\phi} & E & F \\
    \frac{1}{2}G_{\theta} & F & G
  \end{vmatrix}
}{(EG - F^2)^2}
$$ \\
This is so remarkable that Carl Gauss named the result \textit{Theorema Egregium} (in Latin, "remarkable theorem").
\begin{figure}[H]
\centering
\includegraphics[width=4cm, height=5cm]{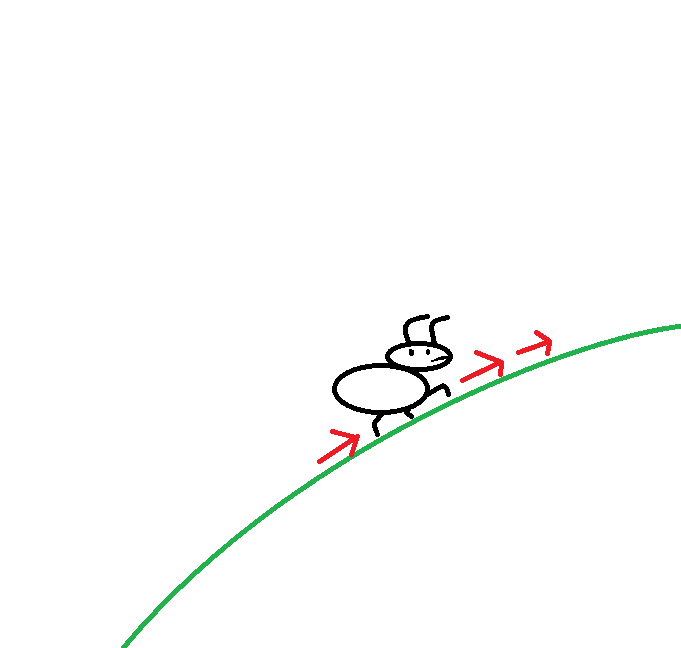}
    \begin{center}
        Figure 2.2: \textit{Theorema Egregium} tells us that simply from the local information gathered by an ant as it traverses the surface, we can determine the Gaussian curvature of a surface.
    \end{center}
\end{figure}
But what about the curvature of an arbitrary smooth (differentiable) manifold with metric $g$, not just surfaces in 3-dimensional Euclidean space? We will generalize the idea of Gaussian curvature to that of \textit{sectional curvature}. Before we can compute the sectional curvature, we must introduce a few key definitions. \\\\
\textbf{Definition 2.5:} atlas \\
A \textit{chart} $(U, \phi)$ on manifold $\mathscr{M}$ consists of an open subset $U$ of $\mathscr{M}$, and a continuous bijection (homeomorphism) $\phi$ from $U$ to an open subset of $\mathbb{R}^n$. An \textbf{atlas} is a collection of charts covering $\mathscr{M}$. \\\\ 
Much like the atlases of our Earth that you may be familiar with, this atlas provides a full coordinate representation of our manifold. A special kind of atlas called a \textit{smooth atlas} does so in a way that allows us to take derivatives on a manifold. \\\\
\textbf{Definition 2.6:} geodesic \\
The shortest distance (locally) between 2 points on a Riemannian manifold. \\\\
On $\mathbb{R}^2$, a geodesic is a straight line. \\
As mentioned previously, we want to examine manifolds equipped with a defined derivative operation. In the Euclidean space, we defined a derivative in the following manner:
$$
f'(x) = \lim_{t \to x} \frac{f(t) - f(x)}{t - x}
$$
Given a curve $c(t)$ lying on a vector field of a manifold $X$, it may seem natural to define to define a derivative at a point $p$ on the curve in the following way:
$$
\lim_{t \to 0} \frac{X(c(t)) - X(p)}{t}
$$
There is a crucial flaw with this definition, however, in that unlike in Euclidean space, there is no ambient space in which the tangent vectors can all reside. A tangent vector at one point on the manifold may reside in a completely different vector space than a tangent vector at another point so that there is no continuity between them. We therefore take the linearization of the manifold at $p$ in the direction of another vector field $Y$:
$$
\nabla_Y X = \lim_{t \to 0} (\frac{X(c(t)) - X(p)}{t})^t
$$
We refer to $\nabla_Y X$ as an \textit{affine connection}. Under our previous definition of the $\nabla$ operator, we could not add 2 tangent vectors along the same curve since they lay in different vector spaces. An affine connection enables us to move vectors from 1 tangent space to another provided that the tangent vectors remain zero with respect to the connection (\textbf{parallel transport}). We have now fully generalized the derivative operation to manifolds. \\\\
There are many different ways to define an affine connection, but the most important for our purposes is the \textit{Levi-Civita connection}: \\\\
\textbf{Definition 2.7:} Levi-Civita connection \\
The unique metric-preserving ($\nabla g = 0$) affine connection that is torsion-free ($\nabla_{X}Y - \nabla_{Y}X - [X,Y] = 0$, where $[X,Y]$ is the \textit{Lie bracket}, $[X,Y] = XY - YX$). \\\\
In keeping with the interpretation of a geodesic being a generalization of straight lines in the Euclidean space to manifolds, we note the following important property of geodesics:
$$\nabla_{\dot{\gamma}} \cdot \dot{\gamma} = 0$$
Using the following properties of affine connections (which should be familiar to students of elementary calculus), we can compute the full meaning of this equation: 
$$ (1) \nabla_X (Y+Z) = \nabla_X Y + \nabla_X Z$$
$$ (2) \nabla_{X+Y} Z = \nabla_X Z + \nabla_Y Z$$
$$ (3) \nabla_{fX} Y = f \nabla_X Y$$
$$ (4) \nabla_X (fY) = f\nabla_X Y + X(f) \cdot Y$$
(For example, (4) is the equivalent of the product rule in calculus). \\\\
In the 2-dimensional case:
$$
   \dot{\gamma} = \Gamma^{1}_{ij} \frac{\partial}{\partial{u}} + \Gamma^{2}_{ij} \frac{\partial}{\partial{v}} 
$$
$$
    \nabla_{\Gamma^{1}_{ij} \frac{\partial}{\partial{u}} + \Gamma^{2}_{ij} \frac{\partial}{\partial{v}}} (\Gamma^{1}_{ij} \frac{\partial}{\partial{u}} + \Gamma^{2}_{ij} \frac{\partial}{\partial{v}}) = 0
$$
By (1):
$$
    \nabla_{\Gamma^{1}_{ij} \frac{\partial}{\partial{u}} + \Gamma^{2}_{ij} \frac{\partial}{\partial{v}}} \Gamma^{1}_{ij} \frac{\partial}{\partial{u}} + \nabla_{\Gamma^{1}_{ij} \frac{\partial}{\partial{u}} + \Gamma^{2}_{ij} \frac{\partial}{\partial{v}}} \Gamma^{2}_{ij} \frac{\partial}{\partial{v}} = 0
$$
By (2):
$$
    \nabla_{\Gamma^{1}_{ij} \frac{\partial}{\partial{u}}} \Gamma^{1}_{ij} \frac{\partial}{\partial{u}} + \nabla_{\Gamma^{2}_{ij} \frac{\partial}{\partial{v}}} \Gamma^{1}_{ij} \frac{\partial}{\partial{u}} + \nabla_{\Gamma^{1}_{ij} \frac{\partial}{\partial{u}}} \Gamma^{2}_{ij} \frac{\partial}{\partial{v}} + \nabla_{\Gamma^{2}_{ij} \frac{\partial}{\partial{v}}} \Gamma^{2}_{ij} \frac{\partial}{\partial{v}} = 0
$$
By (3):
$$
    \Gamma^{1}_{ij} \nabla_{\frac{\partial}{\partial{u}}} \Gamma^{1}_{ij} \frac{\partial}{\partial{u}} + \Gamma^{2}_{ij} \nabla_{\frac{\partial}{\partial{v}}} \Gamma^{1}_{ij} \frac{\partial}{\partial{u}} + \Gamma^{1}_{ij} \nabla_{\frac{\partial}{\partial{u}}} \Gamma^{2}_{ij} \frac{\partial}{\partial{v}} + \Gamma^{2}_{ij} \nabla_{\frac{\partial}{\partial{v}}} \Gamma^{2}_{ij} \frac{\partial}{\partial{v}} = 0
$$
By (4):
$$
    (\Gamma^{1}_{ij})^2 \nabla_{\frac{\partial}{\partial{u}}} \frac{\partial}{\partial u} + \frac{\partial}{\partial u}(\Gamma^{1}_{ij}) \cdot \frac{\partial}{\partial u} + \Gamma^{1}_{ij} \Gamma^{2}_{ij} \nabla_{\frac{\partial}{\partial{u}}} \frac{\partial}{\partial u} + \frac{\partial}{\partial v}(\Gamma^{1}_{ij}) \cdot \frac{\partial}{\partial u} + 
$$
$$
    \Gamma^{2}_{ij} \Gamma^{1}_{ij} \nabla_{\frac{\partial}{\partial u}} \frac{\partial}{\partial v} + \frac{\partial}{\partial u}(\Gamma^{2}_{ij}) \cdot \frac{\partial}{\partial v} + (\Gamma^{2}_{ij})^2 \nabla_{\frac{\partial}{\partial v}} \frac{\partial}{\partial v} + \frac{\partial}{\partial v}(\Gamma^{2}_{ij}) \cdot \frac{\partial}{\partial v} = 0
$$
Given the Levi-Civita connection, and a basis for the tangent space $\{ \partial_i, \partial_j \}$, we define the \textit{Christoffel symbols of the second kind} in the following manner:
$$\nabla_{\partial_i} \partial_j = \Gamma^{k}_{ij} \partial_k$$
i.e., that $\Gamma^{k}_{ij}$ times the basis vector $\partial_k$ is equal to the change in $\partial_j$ in the direction of $\partial_i$. Intuitively, the Christoffel symbols tell us how the tangent basis changes from point-to-point. \\\\
\textbf{Definition 2.8}: Riemann curvature tensor
$$R^i_{jkl} = \Gamma^i_{lj,k} - \Gamma^i_{kj,l} + (\Gamma^i_{kp} \Gamma^p_{lj} - \Gamma^i_{lp} \Gamma^p_{kj})$$
where $\Gamma^i_{lj,k}$ is the partial derivative of $\Gamma^i_{lj}$ with respect to $k$, and $p$ ranges over the components of the manifold. \\\\
\textbf{Definition 2.9}: Sectional curvature
$$K = \frac{\langle R(\phi, \theta)\theta, \phi \rangle}{\langle \phi, \phi \rangle \langle \theta, \theta \rangle - \langle \phi, \theta \rangle^2}$$
We now have the prerequisite background necessary to compute the sectional curvature of the torus. \\\\
Recall the parameterization for the torus given in Definition 2.4. The metric tensor $g_{ij}$ is the same as the first fundamental form matrix.
$$ g_{ij} = 
\begin{pmatrix}
(R + r\cos{\theta})^2 & 0 \\
0 & r^2 \\
\end{pmatrix}
$$
$$
g^{ij} = g_{ij}^{-1} = 
\begin{pmatrix}
\frac{1}{(R + r\cos{\theta})^2} & 0 \\
0 & \frac{1}{r^2} \\
\end{pmatrix}
$$
Let's compute the Christoffel symbols first. We can use the following formula:
$$
\Gamma^{k}_{ij} = \frac{1}{2} g^{kl} \left( \frac{\partial g_{jl}}{\partial x^i} + \frac{\partial g_{il}}{\partial x^j} - \frac{\partial g_{ij}}{\partial x^l} \right)
$$
Notice that $\Gamma^{k}_{ij} = \Gamma^{k}_{ji}$ always. This is equivalent to the Levi-Civita connection being torsion-free. \\\\
Let $1 = \phi, 2 = \theta$:
$$\Gamma^{1}_{21} = \Gamma^{1}_{12} = \frac{-r\sin{\theta}}{R + r\cos{\theta}}$$
$$\Gamma^{2}_{11} = \frac{\sin{\theta}(R + r\cos{\theta})}{r}$$
The reader is encouraged to verify that the remaining Christoffel symbols are all zero. \\\\
We get the following partial derivatives of the nonzero Christoffel symbols:
$$\Gamma^{1}_{21,2} = \frac{-r(R\cos{\theta} + r)}{(R + r\cos{\theta})^2}$$
$$\Gamma^{2}_{11,2} = \frac{R\cos{\theta} + r\cos^2{\theta} - r\sin^2{\theta}}{r}$$
We can rewrite Definition 2.8 to make it easier to compute using the Riemann curvature and metric tensor components and the formula for the dot product (Definition 2.2): 
$$
K = \frac{\phi_j \theta_i R^p_{ijk}g_{lp}. \theta_l \phi_k}{\sum_{i=1}^2 \sum_{j=1}^2 g_{ij}(\theta_i\theta_j + \theta_i\phi_j + \phi_i\phi_j)}
$$
where $p$ ranges over the components of the manifold. \\\\
First, we compute the nonzero components of the Riemann curvature tensor using Definition 2.7 and algebraically manipulate the resulting equations to simplify their form for the next step:
$$
R^{1}_{212} = \frac{rR\cos{\theta} + r^2 - r^2\sin^2{\theta}}{(R + r\cos^2{\theta})}
$$
\begin{align*}
&= \frac{rR\cos{\theta} + r^2 - r^2\sin^2{\theta}}{R + r\cos{\theta}} \cdot \frac{1}{R + r\cos{\theta}} \\
&= \frac{rR\cos{\theta} + r^2(1 - \sin^2{\theta})}{R + r\cos{\theta}} \cdot \frac{1}{R + r\cos{\theta}} \\
&= \frac{rR\cos{\theta} + r^2\cos^2{\theta}}{R + r\cos{\theta}} \cdot \frac{1}{R + r\cos{\theta}} \\
&= \frac{r\cos{\theta}(R + r\cos{\theta})}{R + r\cos{\theta}} \cdot \frac{1}{R + r\cos{\theta}} \\
&= \frac{r\cos{\theta}}{R + r\cos{\theta}} \\
\end{align*}
From the skew-symmetry property of the Riemann curvature tensor:
$$
R^{2}_{212} = -R^{2}_{112} = \frac{\cos{\theta}(R + r\cos{\theta})}{r}
$$
Therefore,
$$K = \frac{g_{11} \cdot R^1_{212}}{g_{11} \cdot g_{22}}$$
$$= \frac{(R + r\cos\theta)^2 \cdot (\frac{r\cos\theta}{R + r\cos\theta})}{r^2(R + r\cos\theta)^2}$$
$$= \frac{r\cos\theta(R + r\cos\theta)}{r^2(R + r\cos\theta)}$$
$$= \frac{\cos\theta}{r(R + r\cos\theta)}$$
As can be seen, the sectional curvature is exactly the same as the Gaussian curvature. This is not a coincidence. The motivation behind sectional curvature arises from taking a tangent plane $\Pi$ to a manifold $\mathscr{M}$ at point $p$, then taking all the geodesics through $p$ whose initial tangent vectors lie on the tangent plane. These form a 2-dimensional manifold $S_{\Pi}$, which inherits its Riemannian metric from $\mathscr{M}$. The sectional curvature is simply the Gaussian curvature of $S_{\Pi}$. As such, in the case of a 2-dimensional surface, there is only one two-dimensional tangent plane at each point, and no choice to make as to what $\Pi$ will be. As such, the sectional curvature is simply the Gaussian curvature of the surface at $p$. \\\\ 
There is something else special about the Levi-Civita connection. Under the Levi-Civita connection, the Christoffel symbols for any parameterization of the Euclidean space are all zero, meaning it has curvature zero, i.e., is a \textit{flat manifold}. While this may be intuitively clear, particularly in the case of $\mathbb{R}^2$, it may not be true given a different choice of affine connection! In fact, we can alternatively define the Euclidean space as being the unique space which has zero curvature under the Levi-Civita connection. In the next chapter, we will be concerned with the "flatness" of specific non-Euclidean manifolds under a different choice of connection.

\section{The Fisher Information Metric}
We now have a solid grounding in the tools of Riemannian geometry relevant in application to statistics. But we have not yet explained why it would even make sense to view probability distributions as a geometric structure, or how such a structure could be constructed. \\\\
For some beginning intuition, let us consider the space of 1-dimensional normal distributions, $\mathcal{N}(\mu, \sigma)$. We know the mean $\mu$ can take on any value whereas the standard deviation $\sigma$ can only take on non-negative values. So for every pair of values $(x, y)$ in the upper half-plane $\mathds{H}$, there exists a corresponding parameterized normal distribution $\mathcal{N}(\mu_x, \sigma_y)$. More precisely, there is a one-to-one correspondence between the set of all normal distributions and the upper half-plane $\mathds{H}$. This provides some intuition for why we may view families of probability distributions as possessing a geometric structure, but to rigorously define the normal distribution family as a manifold, we need to define a metric. \\\\
\textbf{Definition 3.1}: Poincare half-plane model \\
The upper half-plane $\{ (x,y) \mid x,y \in \mathbb{R}, y > 0 \}$ imbued with the Poincare metric $ds^2 = \frac{dx^2 + dy^2}{y^2}$.
\begin{figure}[H]
\centering
\includegraphics[width=7cm, height=9cm]{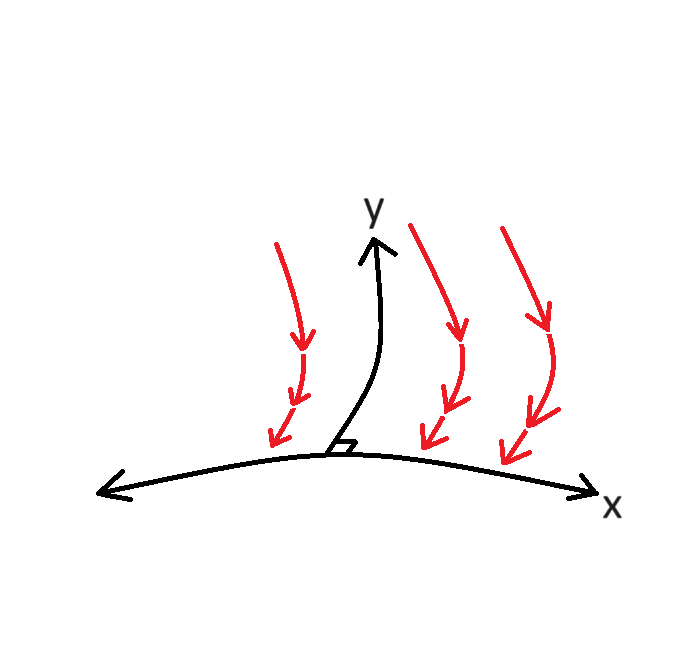}
    \begin{center}
        Figure 3.1: $\mathbb{H}^2$ under the Poincare metric (negative Gaussian curvature)
    \end{center}
\end{figure}

As mentioned in Chapter 2, the upper half-plane imbued with the Euclidean metric has curvature zero. However, surfaces imbued with the Poincare metric have constant negative curvature. This makes the Poincare half-plane a model of $\textit{hyperbolic geometry}$. Objects in the hyperbolic space may not behave the way you are used to. Lines on the half-plane are actually arcs. Distances between points get smaller as you approach the x-axis. If we lived in the hyperbolic plane and tried walking across it in even steps, we would never reach the equator! \\\\
For our purposes, what is important to note is that the half-plane with the Poincare metric is a rigorously defined geometric model. If we can find some natural choice of metric on the set of all normal distributions that is equal to the Poincare metric on the upper half-plane, we can say there is an $\textit{isometry}$ (distance-preserving transformation) between the hyperbolic manifold to the space of normal distributions, i.e., they possess the same geometric structure. \\\\
\textbf{Definition 3.2}: Fisher information \\
The variance of the $\textit{score}$ (the derivative of the log-likelihood (see Definition 1.2) with respect to a parameter $\theta$). Denoted:
$$I(\theta) = \mathbb{E}[(\frac{\partial}{\partial \theta} \log f(X;\theta))^2]$$
(Note: $\mathbb{E}[\frac{\partial}{\partial \theta} \log f(X;\theta)] = 0$) \\\\
Intuitively, Fisher information gives you the expected change in the log-likelihood of a probability distribution with respect to a parameter. If the graph of a distribution is wide and flat (i.e., the derivative of the log-likelihood is small), we can say the random variable $X$ carries little information about the parameter $\theta$. On the other hand, if a distribution has a sharp peak (i.e., the derivative of the log-likelihood is big), we can say that $X$ carries a lot of information about $\theta$, since a change in the value of $\theta$ will clearly result in large changes in the observed values $x$ of $X$. \\\\
\textbf{Definition 3.2.1}: Fisher information matrix \\
The natural extension of Fisher information to the multivariate setting:
$$ [I(\mathbf{\theta})]_{ij} = \mathbb{E}[(\frac{\partial}{\partial \theta_i} \log f(X;\mathbf{\theta}))(\frac{\partial}{\partial \theta_j} \log f(X;\mathbf{\theta}))] $$
Could the Fisher information matrix form the metric tensor that allows us to interpret distributions geometrically? Let us calculate the 2x2 Fisher information matrix for the 1-dimensional normal distribution:
$$f(x; \theta) = \frac{1}{\sigma \sqrt{2\pi}} \exp\left(-\frac{(x - \mu)^2}{2\sigma^2}\right), \mathbf{\theta} = (\mu, \sigma)$$:
$$\log f(x; \theta) = -\frac{1}{2}\log 2\pi  - \frac{1}{2} \log \sigma^2 - \frac{1}{2\sigma^2}(x - \mu)^2$$
$$\frac{\partial}{\partial \mu} \log f(x; \theta) = \frac{x - \mu}{\sigma}$$
$$\frac{\partial}{\partial \sigma} \log f(x; \theta) = \frac{(x - \mu)^2}{\sigma^3} - \frac{1}{\sigma}$$ \\
$$s(\theta) = \begin{pmatrix}
    \frac{\partial}{\partial \mu} \log f(x; \theta) \\\\
    \frac{\partial}{\partial \sigma} \log f(x; \theta)
\end{pmatrix} = 
\begin{pmatrix}
    \frac{x - \mu}{\sigma^2} \\\\
    \frac{(x - \mu)^2}{\sigma^3} - \frac{1}{\sigma}
\end{pmatrix}$$ 
We can also write the Fisher information matrix as the outer product of the score:
$$
I(\mathbf{\theta}) = \mathbb{E}[s(\theta)s(\theta)^T] = \mathbb{E}\begin{pmatrix}
    \frac{(x-\mu)^2}{\sigma^4} & \frac{(x-\mu)^3}{\sigma^5} - \frac{x-\mu}{\sigma^3} \\
    \frac{(x-\mu)^3}{\sigma^5} - \frac{x-\mu}{\sigma^3} & \frac{(x - \mu)^4}{\sigma^6} + \frac{1}{\sigma^2} - \frac{2(x-\mu)^2}{\sigma^4} \\
    \end{pmatrix}
$$
$$
 = 
\begin{pmatrix}
    \mathbb{E}[\frac{(x-\mu)^2}{\sigma^4}] & \mathbb{E}[\frac{(x-\mu)^3}{\sigma^5} - \frac{x-\mu}{\sigma^3}] \\
    \mathbb{E}[\frac{(x-\mu)^3}{\sigma^5} - \frac{x-\mu}{\sigma^3}] & \mathbb{E}[\frac{(x - \mu)^4}{\sigma^6} + \frac{1}{\sigma^2} - \frac{2(x-\mu)^2}{\sigma^4}]
\end{pmatrix}
$$ 
For index position (1,1): 
$$\mathbb{E}[\frac{(x-\mu)^2}{\sigma^4}] = \frac{1}{\sigma^4} \mathbb{E}[(x - \mu)^2] = \frac{1}{\sigma^4} (\mathbb{E}[x^2] -2\mu \mathbb{E}[x] + \mu^2) 
$$
$$
= \frac{1}{\sigma^4}[\mu^2 + \sigma^2 - 2\mu^2 + \mu^2] = \frac{1}{\sigma^4}[\sigma^2] = \frac{1}{\sigma^2}$$ 
For index positions (1,2), (2,1):
$$\mathbb{E}[\frac{(x-\mu)^3}{\sigma^5} - \frac{x-\mu}{\sigma^3}] = \frac{1}{\sigma^5}\mathbb{E}[(x-\mu)^3]$$
$$= \frac{1}{\sigma^5}\mathbb{E}[x^3 -3\mu x^2 +3\mu^2 x + \mu^3] = \frac{1}{\sigma^5}(\mathbb{E}[x^3] - 3\mu\sigma^2 + \mu^3) = \frac{1}{\sigma^5}(\mu^3 + 3\mu \sigma^2 - 3\mu \sigma^2 - \mu^3) = 0$$
For index position (2,2):
$$\mathbb{E}[\frac{(x - \mu)^4}{\sigma^6} + \frac{1}{\sigma^2} - \frac{2(x-\mu)^2}{\sigma^4}] = \frac{1}{\sigma^6}\mathbb{E}[(x - \mu)^4] - \frac{2}{\sigma^4}\mathbb{E}[(x-\mu)^2] + \frac{1}{\sigma^2} = \frac{3}{\sigma^2} - \frac{2}{\sigma^2} + \frac{1}{\sigma^2} = \frac{2}{\sigma^2}$$ \\
\begin{adjustwidth}{-0.5cm}{0cm}
Therefore, the Fisher information matrix for the univariate normal distribution is: \\
\end{adjustwidth}
$$I(\mathbf{\theta}) = \frac{1}{\sigma^2} \begin{pmatrix}
    1 & 0 \\
    0 & 2
\end{pmatrix}$$
This matrix is positive-definite, which satisfies the metric properties. Within a choice of local coordinates, it is exactly equal to the Poincare metric:
$$g_{\text{poincare}} = \frac{1}{\sigma^2}\begin{pmatrix}
    1 & 0 \\
    0 & 1 
\end{pmatrix}$$
This provides both the foundation and the motivation behind information geometry. By applying the Fisher information matrix as a metric (when applicable), we can characterize families of probability distributions as a manifold imbued with this metric tensor. This allows us to consider the informational difference between 2 distributions as a geometric distance. Let us formalize these concepts to extend them beyond the case of the univariate normal distribution. \\\\
\textbf{Definition 3.3}: statistical manifold \\
For parameters $\mathbf{\theta} = (\theta_1, \theta_2, ..., \theta_n)$ lying within a \textit{domain} (open, connected set) $\mathscr{H}$, $\mathscr{M} = \{ p_\theta | \theta \in \mathscr{H} \}$ forms a Riemannian manifold when equipped with the Fisher information matrix as a metric ($\textit{g}_{\textit{FIM}}(\mathbf{\theta})$), provided that: \\
\hspace*{1cm} (1) $\textit{g}_{\textit{FIM}}(\mathbf{\theta})$ is positive-definite for each choice of $\mathbf{\theta} \in \mathscr{H}$. \\
\hspace*{1cm} (2) $\theta \in \mathscr{H} \longmapsto p_{\theta} \in \mathscr{M}$ is a one-to-one map.
\\\\
(1) is necessary for obvious reasons as otherwise the Fisher information matrix (F.I.M.) would not be a metric and the manifold would lack a metric. (2) is necessary because otherwise the likelihood function may offer no information about a given distribution for certain choice of parameters. This means that in our manifold 2 distinct distributions may not have any distance between them, i.e., be identified with the same point. This causes obvious problems. This is known as \textit{non-identifiability}. \\\\
The F.I.M is of supreme importance to information geometry as a result of the following theorem: \\\\
\textbf{Definition 3.4}: sufficiency \\
A statistic $t = T(X)$ is \textit{sufficient} for the underlying parameters of the distribution $\theta$ if the conditional distribution $\mathbb{P}[X|t]$ does not depend on $\theta$. \\\\
In layman's terms: "no other statistic that can be calculated from the same sample provides any additional information as to the value of the parameter" (Fisher, 1922). \\\\
\textbf{Theorem 3.1}: Chentsov's Theorem \\
The Fisher information metric on statistical models
is the only Riemannian metric that is invariant to transformation under sufficient statistics. \\\\
In other words, while other metrics may be computed, the F.I.M is particularly "natural" to the space of statistical manifolds. We will see some powerful applications of this is Chapter 6. \\\\
There is a clever way to rewrite the Fisher information matrix, that will be very useful for our geometric interpretation. First note that: \\
$$\mathbb{E}[\frac{1}{p_\theta} \frac{\partial^2 p_\theta}{\partial\theta_i\partial\theta_j}] = \int \frac{1}{p_\theta} \frac{\partial^2 p_\theta}{\partial\theta_i\partial\theta_j}p_\theta(y) dy = \int \frac{\partial^2 p_\theta}{\partial\theta_i\partial\theta_j} dy$$
$$=\frac{\partial}{\partial\theta_i}\frac{\partial}{\partial\theta_j} \int p_\theta dy = \frac{\partial}{\partial\theta_i}\frac{\partial}{\partial\theta_j} \cdot 1 = 0$$
Now we can rewrite:
$$\mathbb{E}[(\frac{\partial}{\partial \theta_i} \log f(X;\mathbf{\theta}))(\frac{\partial}{\partial \theta_j} \log f(X;\mathbf{\theta}))] = \mathbb{E}[(\frac{1}{p_\theta}\frac{\partial p_\theta}{\partial\theta_i})(\frac{1}{p_\theta}\frac{\partial p_\theta}{\partial\theta_j})] - \mathbb{E}[\frac{1}{p_\theta} \frac{\partial^2 p_\theta}{\partial\theta_i\partial\theta_j}]$$
$$=\mathbb{E}[\frac{1}{p_\theta^2}\frac{\partial p_\theta}{\partial\theta_i}\frac{\partial p_\theta}{\partial\theta_j} - \frac{1}{p_\theta}\frac{\partial^2 p_\theta}{\partial\theta_i\partial\theta_j}] = \mathbb{E}[\frac{-\partial}{\partial\sigma_i}(\frac{1}{p_\theta}\frac{\partial p_\theta}{\partial\theta_j})] = \mathbb{E}[-\frac{\partial^2 \ln p_\theta}{\partial\theta_i\partial\theta_j}]$$
Recall Definition 1.4. Consider the KL divergence between 2 distributions of the same family, $f(x | \theta)$. For any fixed value $\theta'$, the divergence is minimized when $\theta = \theta'$. Now apply the Hessian (second derivative operator) at $\theta = \theta'$ to get: 
$$\nabla_{\theta}^2  \left[ \int p_{\theta'}(x) \log\frac{p_{\theta'}}{p_\theta} dx \right]$$
$$= \int p_{\theta} \cdot \nabla_{\theta}^2 \left[ \log p_{\theta'}(x) - \log p_{\theta}(x) \right] dx$$
$$= -\int p_{\theta} \cdot \nabla_{\theta}^2 \log p_\theta(x) dx$$
We computed above that this is precisely equal to the Fisher information matrix under certain regularity conditions which we will not elaborate on but hold given that $p(y;\theta)$ is a known, specific distribution, which is satisfied in our case. So, each entry in $g_{FIM}$ can be written as the Hessian of the KL divergence. Referring to our knowledge of the second fundamental form from Chapter 2, Fisher information represents the curvature of the KL divergence of a distribution given its parameters. \\\\
Now let's calculate the Hessian of the Euclidean squared distance $D(x,y)^2 = x^2 + y^2$: 
$$\frac{\partial^2}{\partial x^2}D(x,y)^2 = 2$$
$$\frac{\partial}{\partial x}\frac{\partial}{\partial y}D(x,y)^2 = 0$$
$$\frac{\partial^2}{\partial y^2}D(x,y)^2 = 2$$
$$\nabla^2 D(x,y)^2 = 2\begin{pmatrix}
    1 & 0 \\
    0 & 1
\end{pmatrix}$$
So the Hessian of the squared distance produces the metric on the Euclidean space. Like KL divergence, squared distance is non-negative and is zero if and only if $x=y$, and has a positive-definite quadratic form, but does not satisfy the triangle inequality. These measures belong to the class of \textbf{statistical divergences}, but not distance metrics. Similarly, many distances, such as the Euclidean metric, are not valid divergences. The squared distance is the simplest kind of divergence. We have already encountered a member of this class, the Kullback-Leibler divergence from Definition 1.4. Through this lens we may associate the KL divergence on statistical manifolds with the behavior of the squared distance function on the Euclidean space. \\\\
For every statistical divergence $D(\theta': \theta)$, there exists a \textit{reverse divergence} $D^{*}(\theta': \theta) = D(\theta: \theta')$. Divergences occupy an important role in information geometry as a result of the following theorem: \\\\
\textbf{Theorem 3.2}: Divergences induce statistical manifolds \\
$(M, D_g, D_{\nabla}, D^*_{\nabla})$ is an information manifold with: \\
$$^Dg := -\partial_{i,j} D(\theta : \theta_0) \bigg|_{\theta=\theta_0}$$
$$= ^D\Gamma_{ijk} := -\partial_{ij,k} D(\theta : \theta_0) \bigg|_{\theta=\theta_0}$$
$$^{D*}\Gamma_{ijk} := -\partial_{k,ij} D(\theta : \theta_0) \bigg|_{\theta=\theta_0}$$
(The converse holds as well). \\\\
Consider a manifold $\mathscr{M}$. Define 2 global affine coordinate systems $\theta( \cdot )$ and $\eta( \cdot )$ with dual atlases $A = {(\mathscr{M}, \theta)}$ and $A^{*} = {(\mathscr{M}, \eta)}$.
Now consider the following equation $D$: 
$$D(P, Q) = \psi(P) + \phi(Q) - \theta^{i}(P)\eta_{i}(Q)$$
where $\psi, \phi$ are functions such that Hess$^{\nabla}(\psi)$ = Hess$^{\nabla^{*}}(\phi)$ = $g$. The reader is encouraged to verify that $D$ is a valid divergence from the definition. For technical reasons we will not explore further, $D$ defines $\mathscr{M}$ (and vice versa) in a particularly 'nice' way. We call $D$ the \textbf{canonical divergence} of $(\mathscr{M}, g, \nabla, \nabla^{*})$. \\\\
\textbf{Definition 3.5}: Bregman divergence \\
$$B_{F}(\theta : \theta_0) := F(\theta) - F(\theta_0) - (\theta - \theta_0)^\top \nabla F(\theta_0)$$
for a strictly convex function $F$. \\\\
The Bregman divergence can be interpreted as the canonical divergence for a statistical manifold. \\\\
On the statistical manifold induced by a canonical Bregman divergence, $\Gamma_{ijk} = 0, ^{*}\Gamma_{ijk} = 0$, making the manifold both $^{F}\nabla$-flat and $^{F}\nabla^{*}$-flat. This is known as a \textit{dually flat} manifold. As mentioned in Definition 2.5, geodesics are ordinarily the shortest distance between 2 points only locally (consider 2 points on a great circle of a sphere-- the shorter distance between them is a geodesic, but so is the longer distance). On a flat manifold, however, geodesics are the shortest paths between 2 points everywhere, much like lines in the Euclidean space. For a triangle on the manifold formed by 3 geodesic sides, there are 2 choices of geodesic for each side, either $\gamma$ or $\gamma^{*}$, giving us $2^3=8$ geodesic triangles. This leads us to the first major result of this thesis:
\section{The Pythagorean Theorem}
\textbf{Lemma 4.1}: Pythagorean identities 
$$\gamma^{*}(P,Q) \bot_{g} \gamma(Q, R) \Leftrightarrow (\eta(P) - \eta(Q))^T(\theta(Q) - \theta(R)) = (\eta_i(P) - \eta_i(Q))^T(\theta_i(Q) - \theta_i(R)) = 0$$
$$\gamma(P,Q) \bot_{g} \gamma^{*}(Q, R) \Leftrightarrow (\theta(P) - \theta(Q))^T(\eta(Q) - \eta(R)) = (\theta_i(P) - \theta_i(Q))^T(\eta_i(Q) - \eta_i(R)) = 0$$
where $\theta_i, \eta_i$ refer to points along the geodesics $\gamma, \gamma^{*}$ within their respective coordinate systems. An illustration can be seen below. \\\\ 
\textbf{Theorem 4.1}: Pythagorean Theorem of statistical manifolds \\
Let $P,Q,R$ be points on a statistical manifold $\mathscr{M}$ (fixed distributions) such that $PQ$ is a $\nabla$-geodesic and $QR$ is a $\nabla^{*}$-geodesic that form a right angle at $Q$ with respect to $g$ ($\gamma^{*}(P,Q) \bot_{g} \gamma(Q,R)$), then:
$$D(P, R) = D(P, Q) + D(Q, R)$$
where $D$ is the Bregman divergence canonical to $\mathscr{M}$ such that Hess($D$) results in the Fisher information metric. \\\\
\textbf{Proof}: \\
See definition of a canonical divergence. Remember that $D(Q,Q)$ will naturally be 0. For convex $\phi, \psi$: \\
$$
D(P, Q) + D(Q,R) = \psi(P) + \phi(Q) - \theta(P) \cdot \eta(Q) + \psi(Q) + \phi(R) - \theta(Q) \cdot \eta(R)
$$
$$
= \psi(P) + \phi(R) - \theta(P) \cdot \eta(R) + \theta(P) \cdot \eta(R) + \theta(Q) \cdot \eta(Q) - \theta(P) \cdot \eta(Q) - \theta(Q) \cdot \eta(R)
$$
$$
= D(P, R) + (\theta(Q) - \theta(P))^T(\eta(Q) - \eta(R))
$$
The residual term is zero by our assumption. $\blacksquare$ \\\\
So what is the connection between this theorem and the legendary statement of Pythagoras made over 2000 years ago? Recall our identification of the KL divergence on statistical manifolds with the squared Euclidean distance (divergence over Euclidean space). Since the KL divergence is itself a type of Bregman divergence, the theorem is stating that the square of the 'distance' between any point $P$ on a geodesic $\gamma$ and any point $R$ on geodesic $\gamma^{*}$ where $\gamma \bot \gamma^{*}$ equals the square of the distance between $P, Q$ across $\gamma$ plus the square of the distance between $Q, R$ across $\gamma^{*}$. Which is of course nothing more than the Pythagorean Theorem!
\begin{figure}[H]
\centering
\includegraphics[width=4cm, height=4cm]{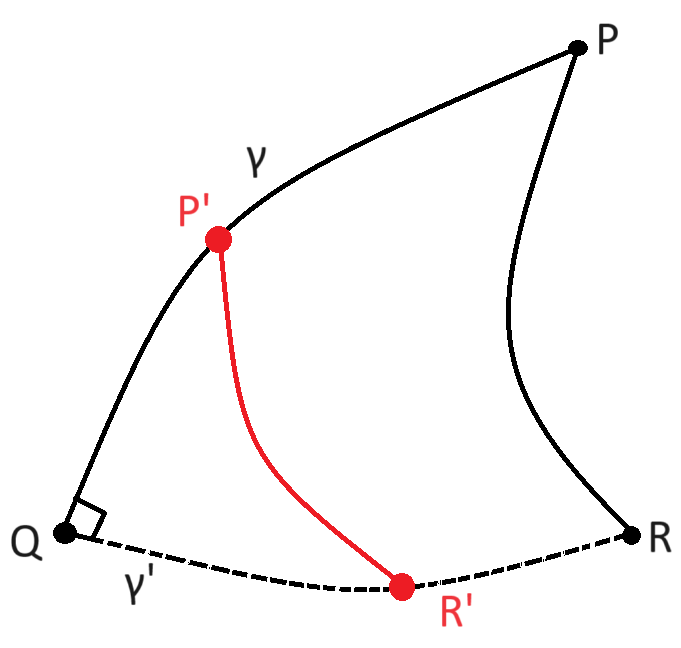}
    \begin{center}
        Figure 4.1: $\gamma(P,Q) \bot_{g} \gamma^{*}(Q, R)$
    \end{center}
\end{figure}
\begin{figure}[H]
\centering
\includegraphics[width=4cm, height=4cm]{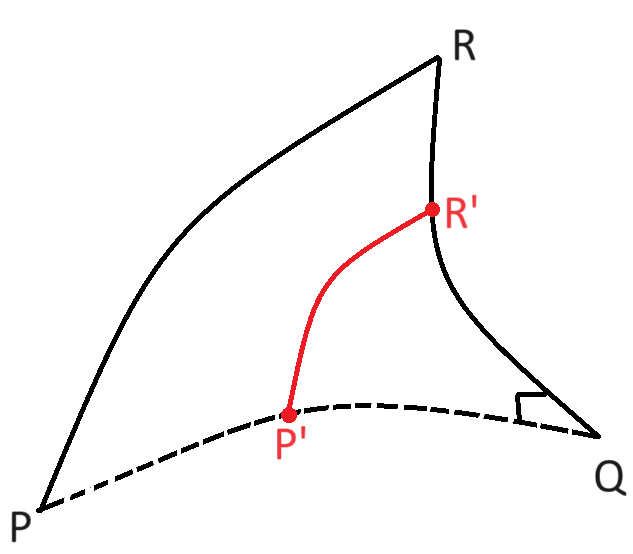}
    \begin{center}
        Figure 4.2: $\gamma^{*}(P,Q) \bot_{g} \gamma(Q, R)$
    \end{center}
\end{figure}

\section{The em Algorithm}
We now return to the setup of our original problem. Given a mix of observed and latent random variables buried within real-world observations, how do we find the correct parameters to fit our model? Let's reframe this scenario within the language of geometry. The set of candidate distributions (distributions of the same family as the ones in our observed data) conditioned by those observed distributions represents our observed data and can be described as detailed above as a manifold $\mathscr{D}$. All points on $\mathscr{D}$ have the same marginal distribution (distributions of known variables). We can write a point $q(\eta)$ on $\mathscr{D}$ as:
$$q(x, z; \eta) = q(x)q(z|x; \eta)$$
where $\eta$ is the unknown parameter of our conditional probability density. 
We wish to identify the points on the manifold of our "true model", $\mathscr{M}$. The intuitive way to do this is to choose the closest point in $\mathscr{M}$ from the data manifold $\mathscr{D}$. We measure the KL divergence between a point $q(\eta)$ on $\mathscr{D}$ and $p(\theta)$ on $\mathscr{M}$:
$$D(q(\eta), p(\theta)) = \int q(x, z; \eta) \log \left( \frac{q(x, z; \eta)}{p(x, z; \theta)} \right) \, dx \, dz$$
We are looking for the $\hat{\theta}, \hat{\eta}$ that minimize $D$. We call this geometric interpretation the \textit{em algorithm} to distinguish it from the original EM algorithm. \\\\
\textbf{Definition 5.1}: $\mathit{m}$-geodesic \\
A set of interior points between 2 distributions $p(x)$ and $q(x)$:
$$r(x;t) = (1 - t) \cdot p(x) + t \cdot q(x), \quad t \in (0, 1)$$
\\\\
\textbf{Definition 5.2}: $\mathit{e}$-geodesic \\
A set of interior points between 2 distributions $p(x)$ and $q(x)$ under a logarithm:
$$\log r(x;t) = (1 - t) \cdot \log p(x) + t \cdot \log q(x) + a(t), \quad t \in (0, 1)$$
where $a(t) = \log \int p(x)^{1-t} \cdot q(x)^t \, dx$ 
\\\\
\textbf{Definition 5.3}: x-projection \\
For point $p$ on manifold $\mathscr{M}$, $\hat{p}$ is the x-projection of $p$ when the $\mathit{x}$-geodesic connecting $p$ to $\hat{p}$ is orthogonal to $\mathscr{M}$ with respect to the FIM $g$ at $\hat{p}$.
\\\\
We may now introduce the \textit{em}-algorithm: \\\\
\textbf{e-step}:\\
Apply the e-projection from $\theta_t$ to $\mathscr{D}$, and obtain $\eta_{t+1}$ by finding 
$$\eta_{t+1} = \arg\max_{\eta} D(q(\eta), p(\theta_t))$$.
\begin{flushleft}
\textbf{m-step}:\\
Apply the m-projection from $\eta_{t+1}$ to $\mathscr{M}$, and obtain $\theta_{t+1}$ by finding 
$$\theta_{t+1} = \arg\max_{\theta} D(q(\eta_{t+1}), p(\theta))$$.
\end{flushleft}
Paralleling the EM algorithm, we iterate until the estimates stop changing.

\begin{figure}[H]
\centering
\includegraphics[width=8cm, height=10cm]{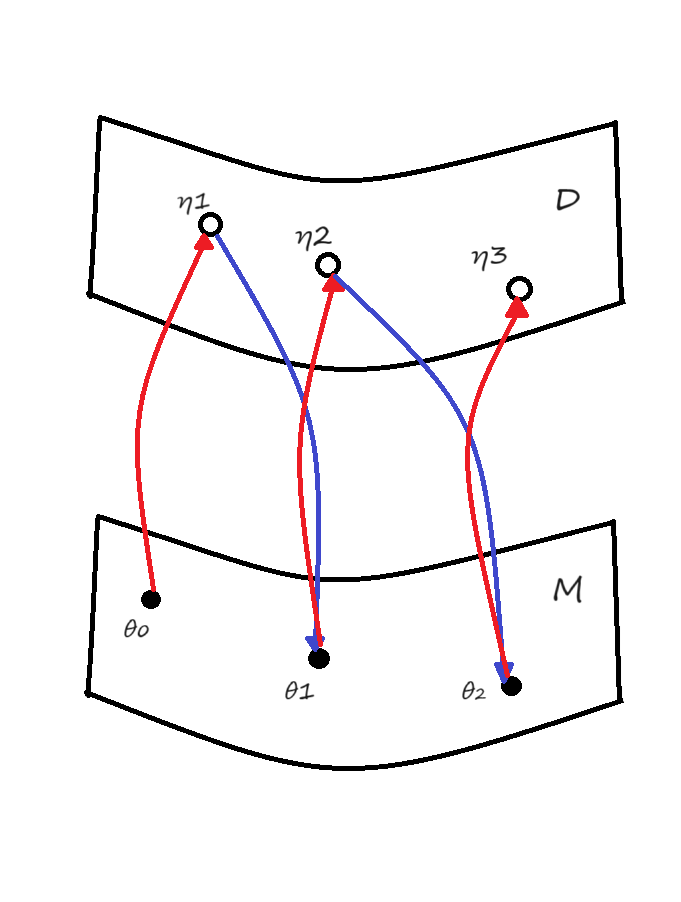}
    \begin{center}
        Figure 5.1: Visualizing $e$- and $m$-projection from $\mathscr{M}$ to $\mathscr{D}$ and vice versa.
    \end{center}
\end{figure}

Comparing Definition 5.1 to Definition 1.3, we see that the $\mathit{e}$-geodesic is an example of a mixture distribution. 
\\\\
We wish to strengthen the connection between the geometric model and the underlying statistics. First, we will introduce a few more key definitions and theorems: \\\\
\textbf{Definition 5.4}: entropy \\
For a discrete probability distribution $P(x)$, the \textit{entropy} is given by: 
$$H(X) = - \sum_{i} P(x_i) \log P(x_i)$$
In the continuous case:
$$H(X) = - \int_{x} P(x_i) \log P(x_i) dx$$
Entropy represents the expected amount of 'surprise' inherent to a random variable. For an example, consider a televised state lottery. To the 999,999 out of a 1,000,000 viewers who lost the draw, the event contains very little surprise. But to the one lucky winner, his personal result was extremely surprising because the probability of that event $E$ was so low, hence, $\log(\frac{1}{p(E)})$ is high. However, the total entropy is balanced out in the average by the low surprisal of the 999,999 losers. \\\\
\textbf{Theorem 5.1}: principle of maximum entropy \\
The distribution which best represents the current state of knowledge about a system is the one with largest entropy. \\\\
Recall what the EM algorithm attempts to achieve, which is to discover the distribution that best models a system given known data. Logically, the result of the EM algorithm should follow the principle of maximum entropy. For computational reasons, we will focus specifically on members of the \textit{exponential family of distributions}, which contains all the most commonly used distributions in applied fields such as ML or statistical physics (normal, binomial, exponental, etc.). \\\\
\textbf{Definition 5.5}: exponential family \\
A member of the exponential family can be expressed in the following form:
$$f(X;\theta) = h(X) \exp\left( \eta(\theta) \cdot T(X) - A(\theta) \right)$$
where $h(X)$ is a non-negative function, $T(X)$ is a sufficient statistic, $\eta(\theta)$ is a function of the parameter, and $A(\theta)$ ensures integration / summation to 1. \\\\
Without data, it is clear that the distribution that maximizes $H(X)$ is the uniform distribution. Given prior information about the distribution, as in our problem setup, however, we can channel this information into one or more constraints, and maximize the entropy over the constraints. Recall that the method of \textbf{Lagrangian multipliers} is used to optimize over a constraint in calculus. \\\\
\textbf{Theorem 5.2}: method of Lagrangian multipliers \\
Solve the following system of equations:
$$\nabla f(x,y,z) = \lambda \cdot \nabla g(x,y,z)$$
$$g(x,y,z) = k \text{    (constraint)}$$
Extend to $n$ many constraint equations. \\\\ 
All of this should provide the necessary background to calculate the $m$-projection step via the below theorem, as may have already become clear. \\\\
\textbf{Theorem 5.3}: \\
To minimize the KL divergence traversed along the $m$-geodesic is equivalent to maximizing the entropy of the points of $\mathscr{M}$.\\\\
\textbf{Proof}: \\
Our constraint equations are:
$$\mathbb{E}[f_k(x)] = \sum^n_{i=1}f_k(x_i)p_i = g_k \in \mathbb{R}$$
$$\sum^n_{i=1}p_i = 1$$
By Theorem 5.2, $\nabla \mathcal{L} = 0$, where:
$$\mathcal{L}(p_1, ..., p_n, \lambda_1, ..., \lambda_n, \mu) = $$
$$-\sum^n_{i=1}p_i\log p_i - \sum^m_{k=1} [\lambda_k \sum^n_{i=1} f_k(x_i)p_i - g_k] - \lambda_0(\sum^n_{i=1}p_i - 1)$$
$$0 = \nabla \mathcal{L} = \sum^n_{i=1} (-\log p_i - 1 - \sum^m_{i=1} \lambda_k f_k(x_i) - (\lambda_0 - 1)) \cdot dp_i$$
(We let $\mu = \lambda_0 + 1$ above). Since we could have optimized with respect to any individual $p_i$: 
$$0 = -\log p_i - 1 - \sum^m_{i=1} \lambda_k f_k(x_i) - (\lambda_0 - 1)$$
$$\log p_i = - \sum^m_{i=1} \lambda_k f_k(x_i) - \lambda_0$$
$$p_i = e^{- \sum^m_{i=1} \lambda_k f_k(x_i)}e^{-\lambda_0}$$
$$1 = \sum^n_{i=1} p_i = \sum^n_{i=1} e^{-\sum^m_{k=1} \lambda_k f_k(x_i)}e^{-\lambda_0}$$
$$Z(\lambda_1, ..., \lambda_m) = e^{\lambda_0} = \sum^n_{i=1} p_i = \sum^n_{i=1} e^{-\sum^m_{k=1} \lambda_k f_k(x_i)}$$
We call $Z$ the \textit{partition function}.
$$p_i = \frac{1}{Z(\lambda_1, ..., \lambda_m)}e^{- \sum^m_{k=1} \lambda_k f_k(x_i)}$$
We can rewrite our constraints $g_k, 1 \leq k \leq m$: 
$$g_k = \sum^n_{i=1} f_k(x_i)p_i = \frac{1}{Z(\lambda_1, ..., \lambda_m)} \sum^n_{i=1} f_k(x_i) Z(\lambda_1, ..., \lambda_m)$$
$$= \frac{-1}{Z(\lambda_1, ..., \lambda_m)} \frac{\partial Z(\lambda_1, ..., \lambda_m)}{\partial \lambda_k}$$
$$= \frac{- \log  Z(\lambda_1, ..., \lambda_m)}{\partial \lambda_k}$$
Notice that the maximized entropy $H(p_1, ..., p_n) = -\sum^n_{i=1} p_i \log p_i$ can be rewritten as a function of $g_1, ..., g_m$, $S(g_1, ..., g_m)$:
$$-\sum^n_{i=1} p_i \log p_i = -\sum^n_{i=1} p_i (\log \frac{1}{Z(\lambda_1, ..., \lambda_m)} - \sum^m_{k=1} \lambda_k f_k(x_i))$$
$$= \sum^n_{i=1} p_i \cdot \log Z(\lambda_1, ..., \lambda_n) + \sum^n_{i=1} [p_i \sum^n_{k=1} \lambda_k f_k(x_i)]$$
$$= \sum^n_{i=1} p_i [\lambda_0 + \sum^n_{i=1} p_i[\sum^m_{k=1} \lambda_k f_k(x_i)]]$$
$$= \lambda_0 + \sum^m_{k=1} \lambda_k \sum^n_{i=1} f_k(x_i)p_i$$ (via linearity of summations)
$$= \log Z(\lambda_1, ..., \lambda_m) + \sum^m_{k=1} \lambda_k g_k$$
Therefore, $S(g_1, ..., g_m) = \log Z(\lambda_1, ..., \lambda_m) + \sum^m_{k=1} \lambda_k g_k$. \\
Take derivative of $S$ with respect to $g_k$:
$$\frac{\partial S}{\partial g_k} = \sum^m_{l=1} \frac{\partial \log Z}{\partial g_k} + \sum^m_{l=1} \frac{\partial \lambda_l g_l}{\partial g_k} + \lambda_k$$
as when $l=k$, $\frac{\partial \lambda_k g_k}{\partial g_k} = \lambda_k$, so we can factor the $k$-th term out of the second summation. By the chain rule we can break the first summation up as follows:
$$= \sum^m_{l=1} \frac{\partial \log Z}{\partial \lambda_l} \cdot \frac{\partial \lambda_l}{\partial g_k} + \sum^m_{l=1} \frac{\partial \lambda_l}{\partial g_k}g_l + \lambda_k$$
$$\frac{\partial \log Z}{\partial \lambda_l} = \frac{1}{Z(\lambda_1, ..., \lambda_m)} \sum^n_{i=1} e^{-\sum^m_{k=1} \lambda_k f_k(x_i)} \cdot -f_l(x_i) = -g_l$$
Therefore,
$$\frac{\partial S}{\partial g_k} = \sum^m_{l=1} -g_l \cdot \frac{\partial \lambda_l}{\partial g_k} + \sum^m_{l=1} g_l \cdot \frac{\partial \lambda_l}{\partial g_k} + \lambda_k$$
$$\frac{\partial S}{\partial g_k} = \lambda_k$$
Recall Definition 1.4, of the KL divergence, which is also known as \textit{relative entropy}, as it measures expected excess 'surprise' from model $Q$ if we assume the actual distribution is $P$. As mentioned previously, without prior information the distribution that maximizes entropy is the uniform distribution. Therefore, we let the second argument to the divergence (representing the unknown statistical manifold $\mathscr{M}$) be the uniform distribution of size $N$:
$$KL(p, U) = \int \log \frac{p(x)}{\frac{1}{N}} \cdot p(x) dx$$
$$= \log N + \int [\log p(x)]p(x) \cdot dx$$
$$= \log N - H(p)$$
for each $p \in \mathscr{D}$. \\\\
It follows that:
\begin{center}
argmin$_{p \in \mathscr{D}} KL(p, U)$ = argmax$_{p \in \mathscr{D}} H(p_1, ..., p_n)$. $\blacksquare$ 
\end{center} 
$$$$
Let $\mathbf{q} = p_1, ..., p_n$. $\mathbf{q(x)} = Ce^{-\sum^m_{k=1} \lambda_k f_k(x_i)}$ for some constant $C$.
$$\mathbf{q(x)} \approx e^{-\sum^m_{k=1} \lambda_k f_k(x_i)}$$
This is an example of the \textit{Boltzmann distribution}, which is a member of the exponential family.
$$\mathcal{L}(q, \lambda_1, ..., \lambda_m, \mu) \approx H(q)$$
$$= \mathbb{E}[-\log q(x_i)] = -\sum^m_{i=1} q(x_i) \cdot \log q(x_i)$$
$$= \sum^m_{i=1} [e^{-\sum^m_{k=1} \lambda_k f_k(x_i)} \sum^m_{k=1} \lambda_k f_k(x_i)]$$
Since $\mathbf{q}$ is a valid discrete distribution function, it sums to 1, therefore:
$$= \sum^m_{k=1} \lambda_k f_k(x_i)$$
Meanwhile: 
$$\frac{1}{m} \sum^m_{i=1} \log q(x_i) = \frac{1}{m} \sum^m_{i=1} \sum^m_{k=1} \lambda_k f_k(x_i) = \sum^m_{k=1} \lambda_k f_k(x_i)$$
So, $\mathcal{L}(q, \lambda_1, ..., \lambda_m, \mu) \approx \frac{1}{m} \sum^m_{i=1} \log q(x_i)$. \\\\
Therefore, we can also say:
\begin{center}
argmax$_{q \in \mathscr{D}} H(q)$ = argmax$_{q \in Q} \sum^m_{i=1} \log q(x_i)$
\end{center} 
\begin{adjustwidth}{5cm}{0cm}
     = argmin$_{q \in Q} \sum^m_{i=1} \log \frac{1}{q(x_i)}$ \\
     = argmin$_{q \in Q} \mathbb{E}_{\hat{p} \in \mathscr{D}}[\log \frac{1}{q(X)}]$ \\
     = argmin$_{q \in Q} \mathbb{E}_{\hat{p} \in \mathscr{D}}[\log \frac{\hat{p}(X)}{q(X)}] + \mathbb{E}_{\hat{p} \in \mathscr{D}}[\log \frac{1}{\hat{p}(X)}]$ \\
     = argmin$_{q \in Q} KL(\hat{p}, q) + H(\hat{p})$ \\
     = argmin$_{q \in Q} KL(\hat{p}, q)$
\end{adjustwidth}
where $Q$ is the set of Boltzmann distributions. \\\\
We can rewrite the result of Theorem 5.3 as: 
\begin{center}
argmin$_{p \in \mathscr{D}} KL(p, U)$ = argmin$_{q \in Q} KL(\hat{p}, q)$
\end{center} 
$$$$
One interpretation of this theorem is that if the uniform distribution is the model that maximizes entropy in the absence of prior information, in the general case, the Boltzmann distribution is the general entropy-maximizing distribution, or, phrased in geometric terms, minimizing the geodesic distance between p on the data manifold and points in U is the same as minimizing the geodesic between $\hat{p}$ on the data manifold and points on the manifold of Boltzmann distributions. \\\\
We can also prove the following theorem regarding the $e$-step: \\\\
\textbf{Theorem 5.4}: \\
Let $P = (\theta^P_{\nu}, \theta^P_{h})$ be a point on the model manifold $\mathscr{M}$, and $Q^* = (\eta^*_{\nu}, \eta^*_{h})$ be the $e$-projection back to the data manifold $\mathscr{D}$, where $\eta^*_{\nu} = \theta^P_{\nu} = \bar{r}_{\nu}$ represents the known variable, and $\eta^*_{h}, \theta^P_{h}$ represents the hidden variable we are attempting to estimate. This hidden part $\theta^P_{h}$ is \textit{invariant} under the $e$-projection (1), furthermore, the conditional probabilities and expectation of $\theta_{h}$ are the same at both $P$ and $Q^*$ (2), (3):
$$(1) \theta^{P}_{h} = \theta^{P}_{h}$$
$$(2) \mathbb{P}(\theta^P_{h} | \theta^P_{\nu}) = \mathbb{P}(\theta^*_{h} | \theta^*_{\nu})$$
$$(3) \mathbb{E}[\theta^P_{h} | \theta^P_{\nu}] = \mathbb{E}[\theta^*_{h} | \theta^*_{\nu}]$$
\textbf{Proof}: \\
Recall definition of canonical divergence:
$$KL(Q^*||P) = \mathbb{E}_{Q^*}[\log \frac{P(\theta^*_{h})}{P(\theta^*_{h})}]$$
$$= (\theta^{*}_{\nu} - \theta^{P}_{\nu}) \cdot \bar{r}_{\nu} + (\theta^*_h - \theta^P_h) \cdot \eta^*_h - \psi(\theta^{*}) - \psi(\theta^P)$$
By the orthogonality of the $e$-projection (Definition 5.3), $Q^*$ minimizes $KL(Q^*||P)$, therefore, $\frac{\partial K(Q||P)}{\partial \eta_h} = 0$.
$$\frac{\partial K(Q||P)}{\partial \eta_h} = \frac{\partial \theta_{\nu}}{\partial \eta_h} \cdot \eta_{\nu} + \frac{\partial \theta_h}{\partial \eta_h} \cdot \eta_h + [(\theta_h - \theta^P_h) \cdot \frac{\partial \eta_h}{\partial \eta_h}] - [\frac{\psi(\theta_Q)}{\partial \theta_{\nu}} \cdot \frac{\partial \theta_{\nu}}{\partial \eta_h}] - [\frac{\partial \psi(\theta_Q)}{\partial \theta_{\nu}} \cdot \frac{\partial \theta_{\nu}}{\partial \eta_h}]$$
$$= \frac{\partial \theta_{\nu}}{\partial \eta_h} \cdot \bar{r}_{\nu} + \frac{\partial \theta_h}{\partial \eta_h} \cdot \eta_h + (\theta_h - \theta^P_h) - [\frac{\psi(\theta_Q)}{\partial \theta_{\nu}} \cdot \frac{\partial \theta_{\nu}}{\partial \eta_h}] - [\frac{\partial \psi(\theta_Q)}{\partial \theta_{\nu}} \cdot \frac{\partial \theta_{\nu}}{\partial \eta_h}]$$
We know that $\frac{\partial \psi(\theta_Q)}{\partial \theta_{\nu}} = \bar{r}_{\nu}, \frac{\partial \psi(\theta_Q)}{\partial \theta_h} = \eta_h$. Therefore, in the case that $Q = Q^*$ and $\frac{\partial K(Q||P)}{\partial \eta_h} = 0$:
$$0 = \frac{\partial \theta_{\nu}}{\partial \eta_h} \cdot \bar{r}_{\nu} + \frac{\partial \theta_h}{\partial \eta_h} \cdot \eta_h + (\theta^{*}_h - \theta^P_h) - [\frac{\psi(\theta_Q)}{\partial \theta_{\nu}} \cdot \frac{\partial \theta_{\nu}}{\partial \eta_h}] - [\frac{\partial \psi(\theta_Q)}{\partial \theta_{\nu}} \cdot \frac{\partial \theta_{\nu}}{\partial \eta_h}]$$
$$\frac{\partial \theta_{\nu}}{\partial \eta_h} \cdot \bar{r}_{\nu} + \frac{\partial \theta_h}{\partial \eta_h} \cdot \eta_h + (\theta^{*}_h - \theta^P_h) = [\frac{\psi(\theta_Q)}{\partial \theta_{\nu}} \cdot \frac{\partial \theta_{\nu}}{\partial \eta_h}] + [\frac{\partial \psi(\theta_Q)}{\partial \theta_{\nu}} \cdot \frac{\partial \theta_{\nu}}{\partial \eta_h}]$$
$$\frac{\partial \theta_{\nu}}{\partial \eta_h} \cdot \bar{r}_{\nu} + \frac{\partial \theta_h}{\partial \eta_h} \cdot \eta_h + (\theta^{*}_h - \theta^P_h) = \bar{r}_{\nu} \cdot \frac{\partial \theta_{\nu}}{\partial \eta_h} + \eta_h \cdot \frac{\partial \theta_{\nu}}{\partial \eta_h}$$
$$(\theta^{*}_h - \theta^P_h) = 0$$
$$\theta^{*}_h = \theta^P_h. \blacksquare$$
Recall Theorem 1.1 (Bayes' law). For a multivariate distribution:
$$\mathbb{P}(y|x;\theta) = \frac{\mathbb{P}(x,y;\theta)}{\int_{y'} \mathbb{P}(x,y;\theta) dy'}$$
Substitute a member of the exponential family (Definition 5.5) for $\mathbb{P}$ (note sufficient statistic $T$ is now a function of multiple random variables, $h$ assumed to be solely a function of x):
$$= \frac{h(x) \exp(\eta(\theta) \cdot T(x,y) - A(\theta))}{\int_{y'} h(x) \exp(\eta(\theta) \cdot T(x,y') - A(\theta)) dy'}$$
Let $\eta(\theta) = (\theta^P_{\nu}, \theta^P_{\nu}), T(x,y) = (\eta_{\nu}, \eta_h), A(\theta) = \psi$:
$$= \frac{\exp(\theta^P_{\nu} \bar{r}_{\nu} + \theta^P_{\mu} \eta_h - \psi)}{\int \exp(\theta^P_{\nu} \bar{r}_{\nu} + \theta^P_{\mu} \eta_h - \psi) d\eta_h}$$
$$= \exp(\log (\frac{\exp(\theta^P_{\nu} \bar{r}_{\nu} + \theta^P_{\mu} \eta_h - \psi)}{\int \exp(\theta^P_{\nu} \bar{r}_{\nu} + \theta^P_{\mu} \eta_h - \psi) d\eta_h}))$$
$$= \exp((\theta^P_{\nu} \bar{r}_{\nu} + \theta^P_{\mu} \eta_h) - \log(\int \exp(\theta^P_{\nu} \bar{r}_{\nu} + \theta^P_{\mu} \eta_h) d\eta_h))$$
$$= \exp(\theta^P_{\nu} \bar{r}_{\nu} + \theta^P_{\mu} \eta_h - \Tilde{\psi})$$
The right-hand side expression depends on $\theta^P_{h}, \bar{r}_{\nu}$ but not $\theta^P_{\nu}$, hence:
$$\mathbb{P}(\theta^P_{h} | \theta^P_{\nu}) = \mathbb{P}(\theta^*_{h} | \theta^*_{\nu})$$
$$\mathbb{E}[\theta^P_{h} | \theta^P_{\nu}] = \mathbb{E}[\theta^*_{h} | \theta^*_{\nu}]. \blacksquare$$ \\\\
Theorem 5.4 is essential because it means that when we project our maximized-likelihood estimate back to $\mathscr{D}$, our estimation for the model is not being altered. \\\\
Similar to as how we proved the equivalency between the maximum entropy principle and the $m$-step of the $em$ algorithm, we would like to prove a similar relationship between the $e$-step and Bayes' law, following from the description of the EM algorithm given in Chapter 1. This would prove a general equivalency between the $em$ and EM algorithms. As it turns out, the 2 algorithms are \textit{not} the same. While the $m$ step and M step are equivalent as proven above, the corresponding $e$ and E steps differ. Nevertheless, we can identify the cases in which the 2 algorithms are indeed equivalent. \\\\
Define:
$$s_Q(r_h | r_{\nu}) = \mathbb{E}[r_h | r_{\nu}]$$, 
$$F: (\eta_{\nu}, \eta_h) \mapsto (\eta_{\nu}, s_Q(\eta_{\nu}))$$
where $\eta_h$ is the unconditional expectation of $r_h$ (the 'true' value). The E step of the EM algorithm gives the point $\hat{Q} \in \mathscr{D}$ with the conditional expectation $s_Q$ as its estimate for the unknown / latent variables. By Theorem 5.4, the conditional expectation at P is equal to the $e$-projection of $Q$ to $\mathscr{D}$, so the E step can be equivalently described in terms of the $e$ step as obtaining $Q^*$ via the E step then transforming through the function $F$. Therefore, if $F$ is the identity, the EM and $em$ algorithms are identical. We can formalize and extend this result even further: \\\\
\textbf{Theorem 5.5}: \\
The EM and $em$ algorithms are equivalent if $s_Q$ is linear. \\\\
\textbf{Proof}: \\
Let $s_Q(r_{\nu}) = a +Br_{\nu}$, where $a$ is constant, $B$ is constant matrix. By linearity of expectation:
$$\eta_{h} = \mathbb{E}_Q[s_Q(r_{\nu})] = s_Q(\mathbb{E}_Q[r_{\nu}]) = s_Q(\eta_{\nu})$$
So $F$ maps $(\eta_{\nu}, \eta_h) \mapsto (\eta_{\nu}, \eta_{\nu})$, i.e. the identity map. \\\\
Incredibly, we can extend this result to any function $f(r)$ where $r$ takes on values over a finite set $K$. We begin with a binary random variable $Y$ where $Y = 0, 1$. Then any function $f$ on $Y$ can be written:
$$f(Y) = (f(1) - f(0)) \cdot Y + f(0)$$
Therefore, the $em$ and EM algorithms are equivalent for any binary classification task. \\
Let's push this linearization "trick" even further. Consider a vector-valued function of $r$, $k(r) = \{ k_m(r) | m \in K \}$ where $K$ finite as mentioned above, $k_m(r) = \delta(r - m)$. Let $f(r) = \sum f(m) \cdot k_m(r)$. So $f$ is linear in the function $k(r)$. Therefore, the $em$ and EM algorithms are also equivalent for the Gaussian mixture model and in Boltzmann machine-based neural networks, among others.
\section{Further Applications}
An emerging development is the utility of information geometry in solving certain problems in \textbf{deep learning}, the branch of artificial intelligence that uses artificial neural  networks. \\\\ 
An \textbf{artificial neural network (ANN)} is a machine learning model inspired by the functioning of biological nervous systems used for the purposes of classification, prediction, and artificial generation. For a basic "feedforward" neural network, the ANN consists of a network of nodes, called neurons or units, organized in layers. These layers typically include an input layer, one or more "hidden" layers (whose output is not directly seen by the user), and an output layer. The input layer consists of each data point as a node Each neuron receives input from the neurons in the previous layer, constructs a linear combination of the input with some weights and a bias term, then processes them through an activation function and produces an output. The end result of all the combined layers is pushed through an \textit{activation function} to produce an output. For example, a model for image classification would use an activation function that sorted the input image into one of the available classes.
\begin{figure}[H]
\centering
\includegraphics[width=12cm, height=8cm]{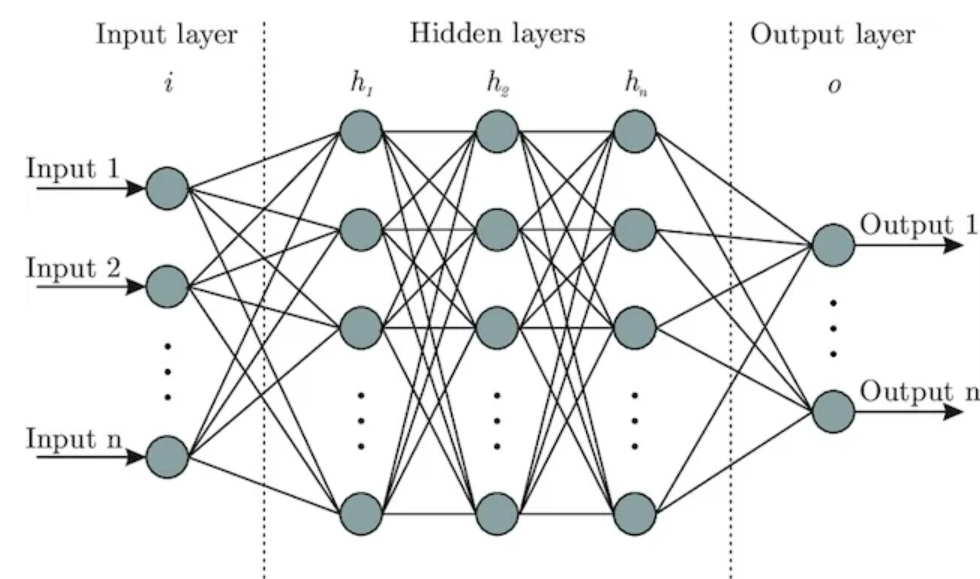}
    \begin{center}
        Feedforward neural network. Credit: u/zyzhu2000, reddit.com/r/MachineLearning
    \end{center}
\end{figure}
\begin{figure}[H]
\centering
\includegraphics[width=7cm, height=8cm]{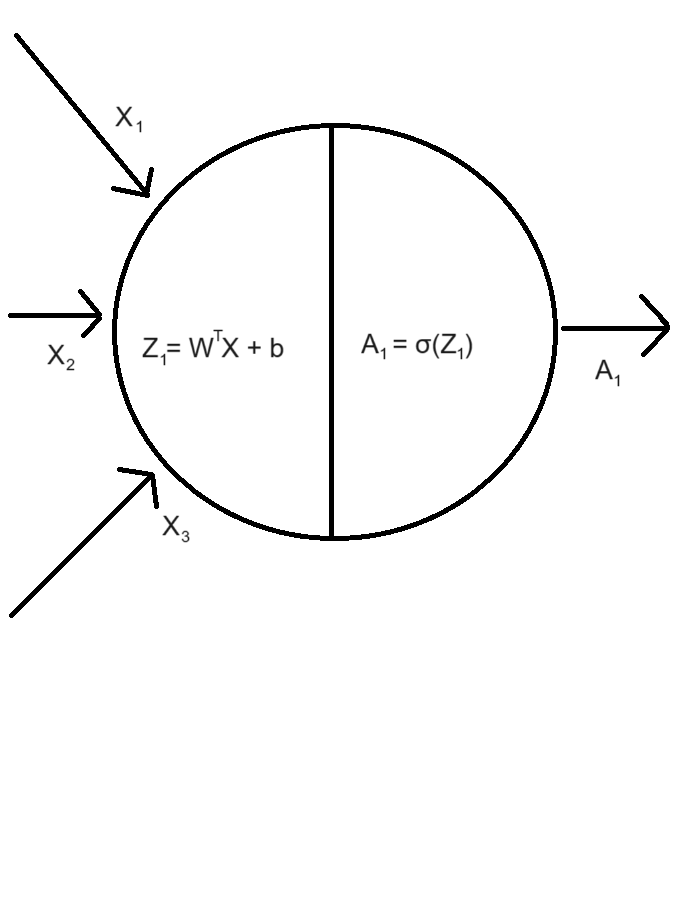}
    \begin{center}
        A single neuron in the above network.
    \end{center}
\end{figure}
To "train" a neural network on a dataset we instantiate the weights of a given architecture of layers. This is done iteratively through the gradient descent algorithm. \\\\
\textbf{Definition 6.1}: method of gradient descent \\
$$W_{t+1}= W_t -\alpha \cdot \nabla \mathcal{L}(W, X)$$
where $W_t$ is the current iteration of our weight vector at step $t$, $\mathcal{L}$ is the \textit{loss function} we are attempting to minimize, X is the entire dataset, and $\alpha$ is some fixed learning rate. We iterate this formula for $n$ steps. \\\\
Intuitively, we are attempting to minimize the loss function by taking $\alpha$-sized steps in the direction of the local minimum. The current industry standard is to only compute the gradient over a small sample of the training data, which is more computationally efficient. This is known as \textit{stochastic gradient descent (SGD)}. \\\\
To obtain the derivatives to apply gradient descent and update the weights and biases, we first take the derivative of the loss function with respect to the output, then progress backwards through the layers using the chain rule. For example:
$$Z_1 = W_1^TX + b$$
$$A_1 = \sigma(Z_1)$$
$$\frac{\partial \mathcal{L}}{\partial W_1} = \frac{\partial Z_1}{\partial W_1} \cdot \frac{\partial A_1}{\partial Z_1} \cdot \frac{\partial \mathcal{L}}{\partial A_1}$$
Gradient descent is a very useful algorithm that has been in use for over a hundred years (Cauchy, 1847). Nevertheless, it comes with several drawbacks. Since gradient descent only searches locally, it can get "trapped" in a local minimum and fail to locate the true minimum of the loss function. Gradient descent is also not stable across different values of $\alpha$. Pick a too-small learning rate and the algorithm may take an infeasibly long time to converge, pick a too-large rate and the algorithm might overstep the true minimum. Suppose the loss function happens to be in terms of Euclidean coordinates, for example, the quadratic loss function:
$$\mathcal{L}(W,X) = c|W_t - W_{\text{optimal}}|^2_2$$
Let $\alpha = \frac{1}{2c}, \frac{\partial \mathcal{L}}{\partial W_t} = 2c|W_t - W_{\text{optimal}}|$, then by Definition 6.1:
$$W_{t+1}= W_t -\frac{1}{2c} \cdot 2c|W_t - W_{\text{optimal}}| = W_{\text{optimal}}$$
The algorithm converges in a single step! Unfortunately, many other loss functions are used in practice which are not Euclidean. A common loss function is the cross-entropy loss for binary classification:
$$\mathcal{L}(A^{[L]}, Y) = -\frac{1}{m} \sum\limits_{i = 1}^{m} (y^{(t)}\log\left(A^{[L] (t)}\right) + (1-y^{(t)})\log\left(1- A^{[L](t)}\right))$$
where $t$ represents the current step of the iteration, $Y$ is the true vector of class labels, and $A^{[L]}$ are the node values of the current layer $L$ of the network. \\
There is no single "correct" choice for $\alpha$ that will guarantee convergence for each weight on every data point because the expression of the loss might have wildly difference magnitudes when evaluated across different $X_i$. In order to solve this problem, we need to stop thinking about minimizing loss in the Euclidean perspective, and iterate instead according to the metric induced by the parameter manifold, that is, the Fisher information matrix. This is as a result of Chentsov's theorem (Theorem 3.1), which indicates that the F.I.M. is the only "natural" metric on a statistical model, so many properties should be describable in terms of this metric. With this insight, we now adapt Definition 6.1 in terms of information geometry: \\\\
\textbf{Definition 6.2}: natural gradient \\
$$W_{t+1} = W_t - \alpha \cdot G^{-1} \nabla \mathcal{L}(W,X)$$
where $G$ is the Fisher information matrix for the parameter vector $W$, defined on the parameter distribution manifold $\mathscr{D}$. \\\\
The natural gradient has many helpful properties that we can exploit: \\\\
\textbf{Definition 6.3}: efficient estimator \\
The minimal variance estimator among all unbiased estimators ($\mathbb{E}[\hat{\theta}] = \theta$), where the lower bound for the variance is given by the \textbf{Cramer-Rao lower bound (CRLB)}, $\text{var}(\hat{\theta}) \geq \frac{1}{I(\theta)}$ \\\\
\textbf{Theorem 6.1}: \\
The natural gradient is asymptotically efficient, that is, as the sample size approaches infinity, its variance-covariance matrix approaches the CRLB. \\\\
However, one major drawback of the natural gradient method is that it requires detailed knowledge on the distribution of our training set to obtain $G^{-1}$. This is a key reason, why as of yet, the natural gradient has yet to make significant headway in industry application, with methods such as the "Adam optimizer" being preferred instead. Adam maintains a different learning rate for each parameter and updates them individually over the course of the training using both the first and second partial derivatives of the loss, as opposed to standard SGD, which is a first-order approximation method, meaning it only uses the first derivative. While the natural gradient explicitly encodes the curvature of the underlying parameter space, Adam implicitly includes this information by updating the parameters differently through the second derivative. Olliver (2016, 2017) has proposed algorithms to improve the process of natural gradient estimation. In 2018, Amari, et al. proved a general analytic form for the inverse Fisher information matrix, eliminating the need for computationally intensive numerical methods. In 2022, Van Sang, et. al developed a new algorithm titled \textit{component-wise natural gradient descent} based on Amari et. al's formula for the inverse FIM, which converged in one-third of the steps as Adam and other industry-standard second-order methods. So while the natural gradient method may only be of theoretical interest now, if further research is done to improve its technical performance, we could see it becoming a genuine competitor to Adam and similar optimization algorithms.
\section{Acknowledgements}
I would like to thank my advisor, Dr. Gunhee Cho, for his kind and patient mentoring of me over the past year. I would also like to thank Malik Tuerkoen and Dr. Saad Mouti for their assistance in helping me understand some particularly challenging concepts in geometry and statistics. Lastly, I would like to thank the friends I've made in the UCSB Mathematics department for their encouragement and camaraderie over the past three years: J.H, G.H, H.Y, and Q.Z-- I definitely would not be graduating without you guys.
\section{References}
(1) Dempster, A. P., Laird, N. M., Rubin, D. B. (1977). Maximum Likelihood from Incomplete Data via the EM Algorithm. Journal of the Royal Statistical Society. Series B (Methodological), Vol. 39, Issue 1, pp. 1-38. \\
(2) Do, Chuong B., Batzoglou, Seraphim (2008). What is the expectation maximization algorithm?, Nature Biotechnology Vol. 26, pp. 897–899 \\
(3) Haugh, Martin. "The EM Algorithm" (2015) URL: \href{http://www.columbia.edu/~mh2078/MachineLearningORFE/EM_Algorithm.pdf}{The EM Algorithm} \\
(4) Shifrin, Theodore. Differential Geometry: A First Course in Curves and Surfaces. 2016. \\
(5) Lee, John M.. Introduction to smooth manifolds. Second. Vol. 218. Graduate Texts in Mathematics. Springer, New York, 2013, pp. xvi+708. ISBN: 978-1-4419-9981-8. \\
(6) Lee, John M.. Riemannian Manifolds: An Introduction to Curvature. Graduate Texts in Mathematics. Vol. 176. New York: Springer-Verlag. 1997. ISBN 978-0-387-98322-6. OCLC 54850593 \\
(7) Schulz, William. "DIFFERENTIAL GEOMETRY ATTACKS THE TORUS" URL: \href{https://www.cefns.nau.edu/~schulz/torus.pdf}{DIFFERENTIAL GEOMETRY ATTACKS THE TORUS} \\
(8) Nielsen, Frank. "An elementary introduction to information geometry". 2018. Retrieved from: arXiv:1808.08271 \\
(9) Cho, Gunhee. "Geometry of 1-dimensional normal distribution with the Fisher information metric" URL: \href{https://drive.google.com/file/d/1PsbBsIjCdCHBP4Hn0Jiegw9wMAmZUmS4/view}{Geometry of 1-dimensional normal distribution with the Fisher information metric} \\
(10) Cho, Gunhee. "KL divergence and Fisher information metric" URL: \href{https://drive.google.com/file/d/1eTd7NUENSAv_ImogjStMQpUfVj_Pa8O6/view}{KL divergence and Fisher information metric} \\
(11) S.Surace "Amari's Pythagorean theorem" URL: \href{https://math.stackexchange.com/questions/3668286/amaris-pythagorean-theorem}{Amari's Pythagorean theorem} \\
(12) Cho, Gunhee. "Maximum entropy principle part 1" URL: \href{https://drive.google.com/file/d/17m7ejn8bh4BhvrPht-61YqphBMqsZUvB/view}{Maximum entropy principle part 1} \\
(13) Cho, Gunhee. "Pythagorean theorem of statistical models" URL: \href{https://drive.google.com/file/d/1sfVnuS2PRwq6GUUNOieVBn65uVvXeZMR/view}{Pythagorean theorem of statistical models} \\
(14) Cho, Gunhee. "Maximum entropy principle (=m projection)" URL: \href{https://drive.google.com/file/d/1Hpk7nY3WqHecKpei0A4CHJ_kHLPhyRIN/view}{Maximum entropy principle (=m projection)} \\
(15) Amari, Shunichi (1994) Information Geometry of the EM and em Algorithms for Neural Networks. Neural Networks, Vol. 8, pp. 1379-1408. \\
(16) AruniRC, "Conditional distribution for Exponential family" URL: \href{https://stats.stackexchange.com/questions/136456/conditional-distribution-for-exponential-family}{Conditional distribution for Exponential family} \\
(17) Amari, Shunichi (1998). Natural Gradient Works Efficiently in Learning. Neural Computation, Col. 10, Issue 2, pp. 251-276. \\
(18) Wild, Cody M. "It's Only Natural: An Excessively Deep Dive Into Natural Gradient Optimization" URL: \href{https://towardsdatascience.com/its-only-natural-an-excessively-deep-dive-into-natural-gradient-optimization-75d464b89dbb}{It's Only Natural: An Excessively Deep Dive Into Natural Gradient Optimization} \\
(19) Amari, Shunichi, Karakida, Ryo, Oizumi, Masafumi. "Fisher Information and Natural Gradient Learning of Random
Deep Networks". 2018. Retrieved from: arXiv:1808.07172 \\
(20) Van Sang, Than, Irvan, Mhd, Yamaguchi, Rie S., Nakata, Toshiyuki, "Component-Wise Natural Gradient Descent - An
Efficient Neural Network Optimization". 2022. Retrieved from: arXiv:2210.05268 \\\\
\textbf{Image Credits}: \\
Piesk, Tilman. "File:Tesseract torus.png" URL: \href{https://upload.wikimedia.org/wikipedia/commons/1/17/Tesseract_torus.png}{File:Tesseract torus.png} \\
u/zyzhu2000. "[D] Best Way to Draw Neural Network Diagrams" URL: \href{https://www.reddit.com/media?url=https%3A%2F%2Fpreview.redd.it%2Fpjj39fhl1pc61.png%3Fwidth%3D640%26format%3Dpng%26auto%3Dwebp%26s%3D37fd1d99b4bf9aa2ff64975ba64daa52a3a141f4}{[D] Best Way to Draw Neural Network Diagrams}
\section{Appendix A: The EM Algorithm code example}
The following code can be accessed and freely used at github.com/SammySuliman. \\\\
We begin by importing some standard libraries for computational mathematics and statistics:
\begin{lstlisting}[language=Python]
    import numpy as np
    import matplotlib.pyplot as plt
    from scipy.stats import norm
    import math
\end{lstlisting}
Initializing a seed for repeatability, we generate normally-distributed GPA data within the range 0.0-4.0 for both populations, then combining the datasets together in a NumPy array. 
\begin{lstlisting}[language=Python]
    np.random.seed(42)
    math_GPAs = np.random.normal(3.7, 0.5, size=20)
    stats_GPAs = np.random.normal(2.8, 0.15, size=20)
    math_GPAs = np.clip(math_GPAs, 0.0, 4.0)
    stats_GPAs = np.clip(stats_GPAs, 0.0, 4.0)
    mixture_data = np.concatenate((math_GPAs, stats_GPAs))
    np.random.shuffle(mixture_data)
\end{lstlisting}
We plot the true distributions of the mixed dataset:
\begin{lstlisting}
    val = 0
    plt.figure(figsize=(8, 6))
    math = plt.plot(math_GPAs, np.zeros_like(math_GPAs) \
    + val, '.', color = 'blue', label='Math Major GPAs')
    stats = plt.plot(stats_GPAs, np.zeros_like(stats_GPAs)\
    + val, '.', color='red', label='Stats Major GPAs')
    plt.legend()
    plt.title('GPA Distribution for Math vs. Statistics Majors')
    plt.xlabel('GPAs')
    plt.show()
\end{lstlisting}
\begin{figure}[H]
\centering
\includegraphics[width=15cm, height=12cm]{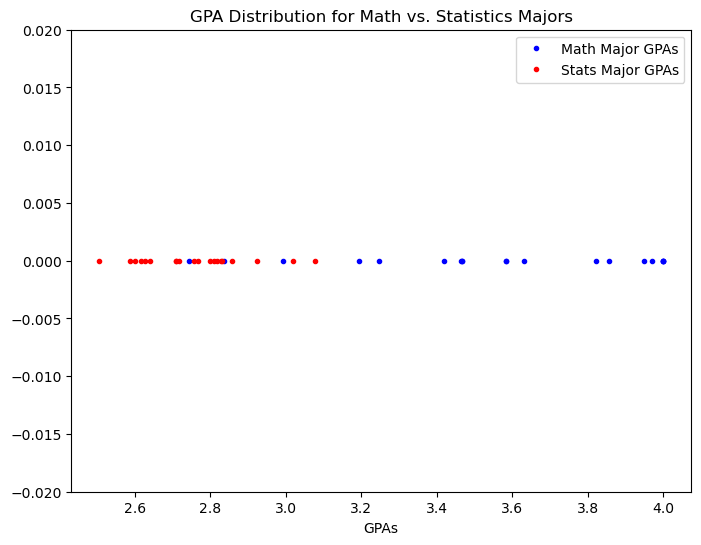}
\end{figure}
As well as our perceived distribution from the point of view of an observer:
\begin{figure}[H]
\centering
\includegraphics[width=15cm, height=12cm]{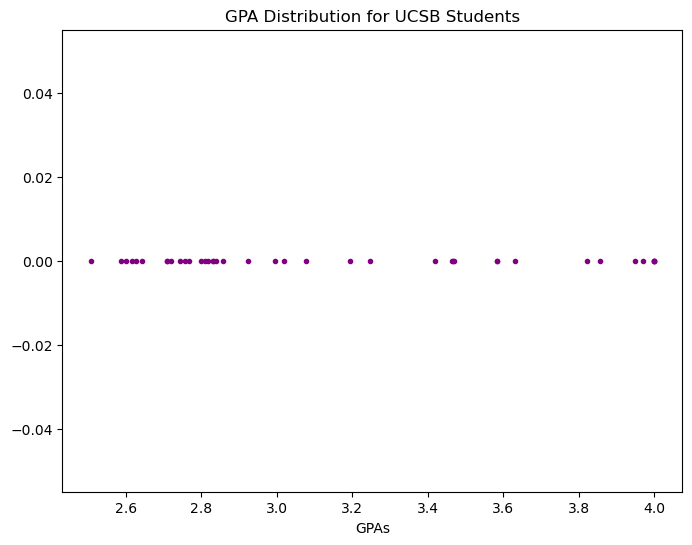}
\end{figure}
We now outline the E-step of the algorithm. The reader is encouraged to compare with the theoretical E-step presented in Chapter 1. The function takes in parameters from the last iteration of the algorithm along with the original dataset. We compute the separate likelihoods under the current guesses of parameters for the mean GPA and stand deviation for math major and statistics major students ($\hat{\mu}_{\text{math}}, \hat{\mu}_{\text{stats}}, \hat{\sigma}_{\text{math}}, \hat{\sigma}_{\text{stats}}$). We then compute $p_{\text{math}}, p_{\text{stats}}$, defined as the proportion of the likelihood for each distribution over the sum of the likelihoods for every individual data point. This represents the known probability that each data point belongs to one distribution or the other at this iteration of the algorithm (the posterior), and return the updated parameters. 
\begin{lstlisting}[language=Python]
    def E_step(data, params):
        (math_likelihood, stats_likelihood, mu_math, mu_stats, \ 
        sigma_math, sigma_stats) = params
        likelihood_math = norm.pdf(data, loc=mu_math, scale=sigma_math)
        likelihood_stats = norm.pdf(data, loc=mu_stats, scale=sigma_stats)
        prob_math = (likelihood_math / (likelihood_math + likelihood_stats))
        prob_stats = (likelihood_stats / (likelihood_math + likelihood_stats))
        return (prob_math, prob_stats, mu_math, mu_stats, \
        sigma_math, sigma_stats)
\end{lstlisting}
For the M-step of the algorithm, we use the analytic solutions for the MLEs of the normal distribution calculated in Chapter 1, progressively weighted on each iteration by $p_{\text{math}}, p_{\text{stats}}$.
\begin{lstlisting}[language=Python]
    def M_step(data, params):
        (prob_math, prob_stats, mu_math, mu_stats, sigma_math,\
        sigma_stats) = params
        mu_hat_math = np.sum(prob_math * data) / np.sum(prob_math)
        sigma_hat_math = np.sqrt(np.sum(prob_math * (data - \
        mu_math)**2) / np.sum(prob_math))
        mu_hat_stats = np.sum(prob_stats * data) / np.sum(prob_stats)
        sigma_hat_stats = np.sqrt(np.sum(prob_stats * (data - \
        mu_stats)**2) / np.sum(prob_stats))
        return (prob_math, prob_stats, mu_hat_math, \
        mu_hat_stats, sigma_hat_math, sigma_hat_stats)
\end{lstlisting}
Putting the complete algorithm together, we randomly initialize our parameters (with set seed) from their expected ranges, then iterate until each of the parameters is within 0.0001 of their previous value or the maximum number of iterations is reached, to prevent a non-terminating loop in the case of stability issues.
\begin{lstlisting}[language=Python]
    def EM_algorithm(data, max_iterations = 100):
        np.random.seed(0)
        mu_math = np.random.uniform(0,4)
        sigma_math = np.random.uniform(0,1)
        mu_stats = np.random.uniform(0,4)
        sigma_stats = np.random.uniform(0,1)
        print(mu_math, sigma_math, mu_stats, sigma_stats)
        num_iterations = 0
        params = (data, data, mu_math, mu_stats, sigma_math, sigma_stats)
        while num_iterations < 1000:
            (prob_math, prob_stats, mu_math, mu_stats, sigma_math, sigma_stats) = params
            params = E_step(data, params)
            params = M_step(data, params)
            num_iterations += 1
            if (math.fabs(params[2] - mu_math) < 0.0001 and \
            math.fabs(params[3] - mu_stats) < 0.0001 \
            and math.fabs(params[4] - sigma_math) < 0.0001 and\
            math.fabs(params[5] - sigma_stats) < 0.0001):
                break
        print(num_iterations)
        return params
\end{lstlisting}
Our results:
$$\hat{\mu}_{\text{math}} = 3.62348312056804$$
$$\hat{\mu}_{\text{stats}} = 2.755661038565841$$
$$\hat{\sigma}_{\text{math}} = 0.3287611585104422$$
$$\hat{\sigma}_{\text{stats}} = 0.12819709371241142$$
Recall the actual values of our parameters:
$${\mu}_{\text{math}} = 3.7$$
$${\mu}_{\text{stats}} = 2.8$$
$${\sigma}_{\text{math}} = 0.5$$
$${\sigma}_{\text{stats}} = 0.15$$
All our estimates are within 0.2 of the actual parameter values, so our algorithm converged successfully.
\section{Appendix B: The natural gradient code example}
The following code example was adapted from Algorithm 1 as presented in \textit{Component-Wise Natural Gradient Descent - An Efficient Neural Network Optimization} by Van Sang et. al. \\\\
To summarize, the paper introduces a version of the natural gradient descent algorithm for convolutional neural networks (neural networks containing layer(s) of filters that are multiplied element-wise by the output). To perform NGD efficiently for a deep network, and avoid inverting a massive Fisher matrix, we use the block-diagonal property of the Fisher information matrix of a neural net, where the F.I.M. for each layer of the network can be represented as a block on the diagonal of the overall F.I.M, and the off-diagonal entries are all zero. From there to invert the F.I.M. for the entire network, we only need to invert the block diagonal entries, which is far more computationally efficient. \\\\
Some explanations of how the variables in the code are defined from the paper:
$$D_a = \text{Cost'}(A_{[L]}, Y) = \frac{\partial \mathcal{L}}{\partial W_{L}} \cdot \frac{\partial W_{L}}{\partial A_{[L]}}$$
$$\frac{\partial A_{[L]}}{\partial W_{L}} = \partial \sigma \cdot A_{[L-1]}$$
where $W_l$ is the weight of the l-th layer, $A_{[l]}$ is the activation of the l-th layer. If we take the activation function to be ReLU (Rectified Linear Unit), the derivative of the activation, $\partial \sigma$ is a Heaviside function.
\begin{lstlisting}[language=Python]
import math
import random
import numpy as np
import torch

class ComponentWiseNaturalGradientDescent:
    def __init__(self, params, activations, layers_dict, \
    lr=0.05, gamma = 0.1, bias=True):
        self.defaults = dict(lr=lr, gamma=gamma, bias=bias)
        self.params = params
        self.layers_dict = layers_dict
        self.activations = activations
        
        # We split the parameters into 2 groups depending on # whether they are from a dense or convolutional
        # layer, so we can update them separately
        
        self.dense_params = []
        self.conv_params = []
        self.conv_activations = []
        self.dense_activations = []
        self.dense_gradients = []
        self.conv_gradients = []

        for p in self.params:
            if p.dim() > 2:
                self.conv_params.append(p)
            else:
                self.dense_params.append(p)

    def zero_grad(self):
    
    # Defining the method to zero out our gradients after
    # each step so they don't accumulate
    
        for param in self.params:
            if param.grad is not None:
                param.grad.detach_()
                param.grad.zero_()

    def step(self, gradients):
        for p in gradients:
            if p.dim() > 2:
                self.conv_gradients.append(p)
            else:
                self.dense_gradients.append(p)

        for a in self.activations:
            if a.dim() > 2:
                self.conv_activations.append(a)
            else:
                self.dense_activations.append(a)

        layer_grad = gradients[-2]
        layer_weight = self.params[-2]
        layer_bias = self.params[-1]
        W_prev = self.activations[-1]
        A_prev = torch.relu(W_prev)
        d_act_L = torch.where((torch.matmul(layer_weight, A_prev.squeeze().T) + layer_bias.unsqueeze(1)) > 0,
                              torch.tensor(1.0), torch.tensor(0.0))
        D_a = torch.matmul(layer_grad, torch.matmul(d_act_L, A_prev.squeeze()).T)
        l2 = len(self.dense_params) // 2
        l3 = len(self.conv_params)
        for l in range(len(list(self.layers_dict.keys())), 1, -1):
            current_layer = list(self.layers_dict.keys())[l-1]
            if current_layer[0:2] == 'fc':
            # Update the dense weights
                layer_grad = self.dense_gradients[2*l2-2]
                layer_weight = self.dense_params[2*l2-2]
                layer_bias = self.dense_params[2*l2-1]
                A_prev = self.dense_activations[l2-1]
                d_act_l = torch.where((torch.matmul(layer_weight, A_prev.T) + layer_bias.unsqueeze(1)) > 0,
                                torch.tensor(1.0), torch.tensor(0.0))
                D_s = torch.matmul(D_a, d_act_l)
                D_a_grad = torch.empty_like(D_s)
                D_a_grad = torch.autograd.grad(outputs=D_s, inputs=D_a, grad_outputs=D_a_grad)
                D_w = torch.matmul(A_prev.T, D_s.T)
                F = torch.matmul(D_w.T, D_w)
                D_w = D_w.unsqueeze(1)
                U = torch.matmul(D_w, torch.inverse(F + self.defaults['gamma'])).squeeze().T
                self.dense_params[2*l2-2].data -= self.defaults['lr'] * U
                l2 -= 2
            if current_layer[0:4] == 'conv':
            # Update the convolutional weights
                D_w = self.conv_gradients[l3-1]
                D_w = D_w.view(D_w.shape[3], -1)
                F = torch.matmul(D_w.T, D_w)
                U = torch.matmul(D_w, torch.inverse(F + self.defaults['gamma']))
                self.conv_params[l3-1].data -= self.defaults['lr'] * U
                l3 -= 1
\end{lstlisting}

Here is the example lightweight 3-layer CNN that we will practice the CW-NGD optimization algorithm on. We use PyTorch as our deep learning software of choice, and train on the classic FashionMNIST dataset for computer vision applications: \\
\begin{lstlisting}[language=Python]
from datetime import datetime
import torch
import torch.nn as nn
import torchvision
import torchvision.transforms as transforms
import torch.optim as optim
from torch.utils.tensorboard import SummaryWriter
from torchviz import make_dot
import os

from ComponentWiseNaturalGradientDescent import  ComponentWiseNaturalGradientDescent

def test_NaturalGradientDescent():
    # Define transforms to preprocess the data
    transform = transforms.Compose(
        [
            transforms.ToTensor(),
        ]
    )

    # Create datasets for training & validation, download if necessary
    training_set = torchvision.datasets.FashionMNIST(
        "./data", train=True, transform=transform, download=True
    )
    validation_set = torchvision.datasets.FashionMNIST(
        "./data", train=False, transform=transform, download=True
    )

    # Create data loaders for our datasets; shuffle for 
    # training, not for validation
    training_loader = torch.utils.data.DataLoader(
        training_set, batch_size=64, shuffle=True
    )
    validation_loader = torch.utils.data.DataLoader(
        validation_set, batch_size=64, shuffle=False
    )

    # Define a simple neural network model
    class SimpleModel(nn.Module):
        def __init__(self):
            super().__init__()
            self.fc1 = nn.Linear(784, 150, bias=True)
            # FashionMNIST is comprised of grayscale images 
            # so the number of input/output channels for the 
            # is 1. We choose a 3x3 kernel. 
            self.conv1 = nn.Conv2d(1, 1, kernel_size=3, \
            stride=1, padding=1, bias=False)
            self.fc2 = nn.Linear(150, 10, bias=True)
            image, _ = training_set[0]
            self.image_size = image.size()

        def forward(self, x):
            x = self.fc1(x)
            x = torch.relu(x)
            x = x.view(-1, self.image_size[0], x.shape[0], x.shape[1])
            x = self.conv1(x)
            x = torch.relu(x)
            x = self.fc2(x)
            x = x.view(x.shape[2], x.shape[3])
            return x
        
    # Define a list to store the linear outputs of the first layer
    preactivations = []

    # Create an instance of the 2-layer model
    device = torch.device("cuda" if torch.cuda.is_available() else "cpu")
    model = SimpleModel().to(device)

    def hook(module, input, output):
        # Save output of the layer
        preactivations.append(output)

    # Registering hook to layer1, layer2
    hook_handle = model.fc1.register_forward_hook(hook)
    hook_handle2 = model.conv1.register_forward_hook(hook)
    
    # Define loss function and optimizer
    criterion = nn.CrossEntropyLoss()
    parameters = list(model.parameters())

    # Create an OrderedDict to store the layers in order
    layers_dict = torch.nn.ModuleDict()

    # Iterate through the modules and add them to the OrderedDict
    for name, module in model.named_modules():
        if name != "":
            layers_dict[name] = module
    optimizer = \
    ComponentWiseNaturalGradientDescent(parameters, activations, layers_dict, lr=0.01)

    running_vloss = 0
    for i, vdata in enumerate(validation_loader):
        vinputs, vlabels = vdata
        vinputs = vinputs.view(vinputs.size(0), -1)
        # Flatten the input images
        voutputs = model(vinputs)
        vloss = criterion(voutputs, vlabels)
        running_vloss += vloss

    def train_one_epoch(epoch_index, tb_writer):
        running_loss = 0.0
        last_loss = 0.0

        # Here, we use enumerate(training_loader) instead of
        # iter(training_loader) so that we can track the batch
        # index and do some intra-epoch reporting
        for i, data in enumerate(training_loader):
            # Every data instance is an input + label pair
            inputs, labels = data
            # Flatten the input images
            inputs = inputs.view(inputs.size(0), -1)

            # Zero your gradients for every batch!
            optimizer.zero_grad()

            # Make predictions for this batch
            outputs = model(inputs)

            # Compute the loss and its gradients
            loss = criterion(outputs, labels)
            loss.backward(retain_graph = True)
            # Gradient Norm Clipping to prevent exploding gradients
            torch.nn.utils.clip_grad_norm_(model.parameters(), max_norm=0.5)
            gradients = \
            torch.autograd.grad(loss, model.parameters(), retain_graph=True)
            # Adjust learning weights
            optimizer.step(gradients)
            # Gather data and report
            running_loss += loss.item()
            if i % 937 == 0:
                last_loss = running_loss / 937  # loss per batch
                print("  batch {} loss: {}".format(i + 1, last_loss))
                tb_x = epoch_index * len(training_loader) + i + 1
                tb_writer.add_scalar("Loss/train", last_loss, tb_x)
                running_loss = 0.0

        return last_loss

    # Initializing in a separate cell so we can
    # easily add more epochs to the same run
    timestamp = datetime.now().strftime("%Y%m%d_%H%M%S")
    writer = SummaryWriter("runs/fashion_trainer_{}".format(timestamp))
    epoch_number = 0

    EPOCHS = 2

    best_vloss = 1_000_000.0

    for epoch in range(EPOCHS):
        print("EPOCH {}:".format(epoch_number + 1))

        # Make sure gradient tracking is on,
        # and do a pass over the data
        model.train(True)
        avg_loss = train_one_epoch(epoch_number, writer)

        running_vloss = 0.0
        # Set the model to evaluation mode,
        # disabling dropout and using population
        # statistics for batch normalization.
        model.eval()

        # Disable gradient computation and reduce memory consumption.
        with torch.no_grad():
            for i, vdata in enumerate(validation_loader):
                vinputs, vlabels = vdata
                vinputs = vinputs.view(vinputs.size(0), -1) 
                # Flatten the input images
                voutputs = model(vinputs)
                vloss = criterion(voutputs, vlabels)
                running_vloss += vloss

        avg_vloss = running_vloss / (i + 1)
        print("LOSS train {} valid {}".format(avg_loss, avg_vloss))

        # Log the running loss averaged per batch
        # for both training and validation
        writer.add_scalars(
            "Training vs. Validation Loss",
            {"Training": avg_loss, "Validation": avg_vloss},
            epoch_number + 1,
        )
        writer.flush()

        # Track best performance, and save the model's state
        if avg_vloss < best_vloss:
            best_vloss = avg_vloss
            model_path = "model_{}_{}".format(timestamp, epoch_number)
            torch.save(model.state_dict(), model_path)

        epoch_number += 1

test_NaturalGradientDescent()
\end{lstlisting}

Unfortunately, we did not achieve good loss values with this choice of optimization algorithm. More research is needed to make NGD applicable in real world application.
\end{document}